\documentclass[12pt]{article}
\usepackage{amssymb, amsmath, amsthm}
\usepackage[all]{xy}

\renewcommand{\baselinestretch}{1.2}

\newcommand{\mr}[1]{\mathrm{#1}}
\newcommand{\mf}[1]{\mathfrak{#1}}
\newcommand{\mc}[1]{\mathcal{#1}}

\newcommand{\Z}{{\bf Z}}
\newcommand{\Q}{{\bf Q}}
\newcommand{\zp}{{\bf Z}_p}

\newcommand{\eq}[1]{(\ref{#1})}

\DeclareMathOperator{\Gal}{Gal} \DeclareMathOperator{\Tra}{Tra}
\DeclareMathOperator{\Res}{Res} \DeclareMathOperator{\Inf}{Inf}
\DeclareMathOperator{\Hom}{Hom} 

\def\cl{\mathrm{Cl}}
\def\Scl#1{\cl_{#1,S}}

\newcommand{\OS}[1]{\mathcal{O}_{{#1},S}}
\newcommand{\US}[1]{\mathcal{O}_{{#1},S}^{\times}}

\newcommand{\Norm}{N_{L/K}}
\DeclareMathOperator{\rank}{rank}

\newcommand{\co}{\, , \,}

\newcommand{\brho}{\bar{\rho}}

\newtheorem{theorem}{Theorem}[section]
\newtheorem{proposition}[theorem]{Proposition}
\newtheorem{lemma}[theorem]{Lemma}
\newtheorem{corollary}[theorem]{Corollary}
\newtheorem{conjecture}[theorem]{Conjecture}
\newtheorem*{thm}{Theorem A}

\theoremstyle{definition}

\newtheorem*{remark}{Remark}

\theoremstyle{remark}
\newtheorem*{acknowledgments}{Acknowledgments}

\numberwithin{equation}{section}

\begin{document}
\title{Massey Products and Ideal Class Groups}
\author{Romyar T. Sharifi}
\date{}
\maketitle

\renewcommand{\baselinestretch}{1}
\begin{abstract}
We consider certain Massey products in the cohomology of a Galois
extension of fields with coefficients in $p$-power roots of unity.
We prove formulas for these products both in general and in the
special case that the Galois extension in question is the maximal
extension of a number field unramified outside a set of primes $S$
including those above $p$ and any archimedean places.  We then
consider those ${\bf Z}_p$-Kummer extensions $L_{\infty}$ of the
maximal $p$-cyclotomic extension $K_{\infty}$ of a number field
$K$ that are unramified outside $S$.  We show that Massey products
describe the structure of a certain ``decomposition-free" quotient
of a graded piece of the maximal unramified abelian pro-$p$
extension of $L_{\infty}$ in which all primes above those in $S$
split completely, with the grading arising from the augmentation
filtration on the group ring of the Galois group of
$L_{\infty}/K_{\infty}$.  We explicitly describe examples of the
maximal unramified abelian pro-$p$ extensions of unramified
outside $p$ Kummer extensions of the cyclotomic field of all
$p$-power roots of unity, for irregular primes $p$.
\end{abstract}
\renewcommand{\baselinestretch}{1.3}

\section{Introduction} \label{intro}

\subsection{Iwasawa modules over Kummer extensions}

Classical Iwasawa theory takes place over the cyclotomic $\zp$-extension $K_{\infty}$ of a number field $K$, for a prime $p$.  One considers finitely generated modules for the Iwasawa algebra $\Lambda = \zp[[\Gal(K_{\infty}/K)]]$, known as Iwasawa modules.  This paper is concerned with Iwasawa theory over a $\zp$-extension $L_{\infty}$ of $K_{\infty}$
that is Galois over $K$.  The pro-$p$ group $\mc{G} = \Gal(L_{\infty}/K)$ is in general nonabelian, being isomorphic to a semi-direct product of two copies of $\zp$: its normal subgroup $G = \Gal(L_{\infty}/K_{\infty})$ and any lift of its quotient $\Gamma = \Gal(K_{\infty}/K)$.  
Whereas $\Lambda$ is isomorphic to a power series ring in one variable over $\zp$, the Iwasawa algebra $\zp[[\mc{G}]]$ is isomorphic to a ring of $\zp$-power series in two variables  that do not commute if $\mc{G}$ is nonabelian.  The structure theory of finitely generated $\zp[[\mc{G}]]$-modules is therefore considerably more complicated than in the classical case (see \cite{venjakob} for a study).

Our idea is a simple one.  Given a finitely generated $\zp[[\mc{G}]]$-module $M$, it is possible to create an infinite sequence of $\Lambda$-modules, each a quotient of the last, as follows.  Let $I_G$ denote the augmentation
ideal in $\zp[[G]]$, and consider the quotients 
$$
	\mr{gr}^k_G M  = I_G^k M/I_G^{k+1} M
$$
for $k \ge 0$. 
The graded object
$$
	\mr{gr}_G M = \bigoplus_{k=0}^{\infty} \mr{gr}^k_G M
$$
is a module over the graded ring $\mr{gr}_G \zp[[\mc{G}]]$.
In particular, each $\mr{gr}^k_G M$ is a finitely generated $\Lambda \cong \mr{gr}_G^0 \zp[[\mc{G}]]$-module.  Since
$$
	\mr{gr}^k_G \zp[[G]] = I_G^k/I_G^{k+1} \cong G^{\otimes k},
$$ 
the action of its first graded piece provides a surjection 
$$
	\mr{gr}^k_G M \otimes_{\zp} G \twoheadrightarrow \mr{gr}^{k+1}_G M
$$ 
of $\Lambda$-modules.  (Note that $G$ is isomorphic as a $\Lambda$-module to a twist of $\zp$ by a possibly nonintegral power of the cyclotomic character.)  Of course, a good deal of information on $M$ can be lost in the passage to the graded module, but much is retained.  For instance, we can detect the 
torsionness of $M$ as a $\zp[[G]]$-module from the 
finiteness of the $\mr{gr}^k_G M$ 
for large  $k$.

Now, one of the most typical Iwasawa modules to consider is the Galois group $\mc{Z}_K$ of the maximal unramified abelian pro-$p$ extension of $K_{\infty}$.  The action of $\Gamma$ on $\mc{Z}_K$ is given by lifting and conjugating.  By class field theory, $\mc{Z}_K$ can also thought of as the inverse limit of $p$-parts of class groups of intermediate subfields of $K$, and the action is induced by the usual action of $\Gamma$ on ideals.  By a standard result of Iwasawa, $\mc{Z}_K$ is a finitely generated, torsion $\Lambda$-module.

We wish to consider $\mc{Z}_L$, the Galois group of the maximal unramified abelian pro-$p$ extension of $L_{\infty}$, and its relation to $\mc{Z}_K$.    Let us assume that $L_{\infty}/K$ is unramified outside a finite set of primes of $K$.  We remark that $\mc{Z}_L/I_G\mc{Z}_L$ is closely related to $\mc{Z}_K$, mapping to it with kernel and cokernel finitely generated over $\zp$.  As a consequence, $\mc{Z}_L$ is finitely generated over $\zp[[\mc{G}]]$.  
Assume now that $K_{\infty}$ contains all $p$-power roots of $1$.
We show in this paper that the higher graded quotients $\mr{gr}^k_G \mc{Z}_L$ for positive $k$ relate to quotients of $\mc{Z}_K$ by the image of certain operations in Galois cohomology called Massey products.  To be precise in our description, we must replace $\mc{Z}_K$ and $\mc{Z}_L$ with closely related Iwasawa modules.

We may consider the quotient of $\mc{Z}_K$ that is the Galois group $\mc{A}_K$ of the maximal unramified abelian pro-$p$ extension of $K_{\infty}$ in which all primes above those in $p$ split completely.  We also have the analogously defined Galois group $\mc{A}_L$ over $L_{\infty}$.   
Again by class field theory, the modules $\mc{A}_K$ and $\mc{A}_L$ are isomorphic to the inverse limits of $p$-parts of the quotients of class groups by the classes of primes above  $p$ in the intermediate subfields of $K_{\infty}$ and $L_{\infty}$, respectively.  It is the higher graded quotients of $\mc{A}_L$ that relate to Massey products, as we shall describe below.  We now state our main result (Theorem \ref{mainthm}), with a strong extra restriction for simplicity of presentation.

\begin{thm}
  Assume that there is a unique prime above $p$ in $L_{\infty}$ and that
  $L_{\infty}/K_{\infty}$ is unramified outside of it.  For each $k \ge 1$,
  there is a canonical isomorphism of $\Lambda$-modules,
  $$
    \mr{gr}^k_G \mc{A}_L
    \cong (\mc{A}_{K}/\mc{P}_{L/K}^{(k)})
    \otimes_{\Z_p} G^{\otimes k},
  $$
  where $\mc{P}_{L/K}^{(k)}$ is a group of inverse limits of $(k+1)$-fold
  Massey products with $k$-copies of a Kummer generator of $L_{\infty}$.
\end{thm}

Some applications of this result are given in Section \ref{examples} in the case that $K = \Q(\mu_p)$ with $p$ odd.  When $L_{\infty}$ is defined by a sequence of cyclotomic $p$-units in $K_{\infty} = \Q(\mu_{p^{\infty}})$, we show that $\mr{gr}^1_G \mc{A}_{L}$ may be either zero or nonzero by providing examples of both instances when $p = 37$ (Propositions \ref{p37} and \ref{larger}).  For such $L_{\infty}$, we conjecture that $\mc{A}_L$ is finitely generated as a $\zp[[G]]$-module (Conjecture \ref{unitconj}).   
When $L_{\infty}$ is a CM-field, however, we show that $\mc{A}_L$ is not always finitely generated over $\zp[[G]]$ (Proposition \ref{pseudo}).  
We have good reason to believe that the Massey products of cyclotomic $p$-units will provide highly interesting arithmetic objects with an interpretation in terms of $p$-adic $L$-functions of modular forms, but we leave that for future endeavors.

\subsection{Construction of Massey products}

In order to describe the Massey products of interest, it is easiest if we begin by working in a general context.  For further technical details on the constructions we give in this subsection, see Sections \ref{embed}-\ref{massey}.
Let $H$ be a group, and let $N \ge 1$.  Massey products are products of two or more elements of $H^1(H,\Z/N\Z)$ that take values in a quotient of $H^2(H,\Z/N\Z)$.

The simplest form of Massey product is the cup product.  Let $\chi_1, \chi_2$ be two characters
in $H^1(H,\Z/N\Z)$.  Recall that the cup product $\chi_1 \cup \chi_2 \in
H^2(H,\Z/N\Z)$ is represented by the two-cocycle $\nu \colon (h_1,h_2) \mapsto \chi_1(h_1)\chi_2(h_2)$ for $h_1, h_2 \in H$.  The cup product provides the obstruction
to the existence of a homomorphism $\rho \colon H \to GL_3(\Z/N\Z)$ with
$$
	\rho(h) = \left( \begin{matrix} 1 & \chi_1(h) & \kappa(h) \\
	0 & 1 & \chi_2(h) \\ 0 & 0 & 1 \end{matrix} \right)
$$
for some function $\kappa \colon H \to \Z/N\Z$. 
Put precisely, the two-coboundary of $\kappa$ is $-\nu$, so the existence of $\rho$ is equivalent to $\nu$ being a coboundary.  Note that there is some ambiguity in the choice of $\kappa$, when it exists, being well-defined only up to a character in $H^1(H,\Z/N\Z)$. 

In general, we can define Massey products recursively.  Let $\chi_1, \chi_2,
\ldots, \chi_q$ be characters in $H^1(H,\Z/N\Z)$.  Suppose we have a map
$\rho \colon H \to GL_{q+1}(\Z/N\Z)$
into the group of unipotent upper-triangular matrices $T_{q+1}$, with
\begin{equation} \label{masseymatrix}
  \rho(h) = \left( \begin{array}{cccccc} 1 & \chi_1(h) &
  \kappa_{1,3}(h) & \ldots & \kappa_{1,q+1}(h) \\
  & 1 & \chi_2(h)  & \ldots & \kappa_{2,q+1}(h) \\
   &  & \ddots & \ddots & \vdots \\
  & &  & 1 & \chi_{q}(h) \\
  & &  &  & 1 \end{array} \right)
\end{equation}
and such that ``modulo the upper right-hand corner," $\rho$ is a
homomorphism.  That is to say, letting $Z_{q+1}$ denote the subgroup of $T_{q+1}$ consisting of matrices with with zeros everywhere off of the diagonal except
the upper right-hand corner, we require that the composition $\bar{\rho}$ of $\rho$ with the surjection $T_{q+1} \to T_{q+1}/Z_{q+1}$ is a homomorphism.  Then
\begin{equation} \label{matrixcocycle}
  \nu \colon (h,h') \mapsto \chi_1(h)\kappa_{2,q+1}(h') +
  \kappa_{1,3}(h)\kappa_{3,q+1}(h') + \ldots +
  \kappa_{1,q}(h)\chi_{q}(h'),
\end{equation}
is a $2$-cocycle.  If $\rho$ is a homomorphism, then $\nu$ is the negative of the coboundary of $\kappa_{1,q+1}$.  Conversely, if $\nu$ has trivial class, then there exists a lift of $\bar{\rho}$ to a homomorphism to $T_{q+1}$.

Ideally, the cohomology class of the $2$-cocycle $\nu$ of \eqref{matrixcocycle} would represent the cohomology class of the $q$-fold Massey product of the $\chi_1, \ldots, \chi_q$.  However, the ambiguity in the choice of the $\kappa_{i,j}$ defining $\bar{\rho}$ means that this class is not well-defined.  For instance, in the case that $q = 3$, we could modify 
$\kappa_{1,3}$ by a cocycle $\psi$, and this will modify $\nu$ by the $2$-cocycle representing $\psi \cup \chi_3$.  Similarly, changing $\kappa_{2,4}$ by a cocycle modifies $\nu$ by the cup product of $\chi_1$ with a cocycle.  The Massey triple product is therefore taken to be the class of $\nu$ in the quotient
$$
	H^2(H,\Z/N\Z)/(H^1(H,\Z/N\Z) \cup \chi_3 + \chi_1 \cup H^1(H,\Z/N\Z)).
$$
A similar story holds for higher-fold Massey products.  In general, the
choice of the maps defining $\bar{\rho}$ (the $\kappa_{i,j}$ with $1 \le i < j-1 \le q$ and $(i,j) \neq (1,q+1)$) is called a {\em defining system} of the Massey product of $\chi_1, \ldots \chi_q$.  The indeterminacy in the Massey product arises from the choice of defining system.  

In this paper, we specialize to the case of $q$-fold Massey products in which all the $\chi_i$ are equal, except possibly $\chi_q$.  Let $N = p^m$ for some $m \ge 1$.  Assume that we are given a fixed surjective homomorphism $\chi \colon H \to \Z/p^n\Z$ for some $n$ with $q \le p^{n-m+1}$,
and let $\chi_i$ be the reduction of $\chi$ modulo $p^m$ for $i \le q-1$.   Consider the subgroup $C_{q+1}$ of $T_{q+1}$ consisting of unipotent matrices with zeros outside the rightmost column.  Our assumptions imply that there exist canonical choices of the maps $\kappa_{i,j}$ for $j \le q$ such that $\rho$ composed with $T_{q+1} \to T_{q+1}/C_{q+1}$ is a homomorphism with cyclic image.  
In fact, we take
$$
	\kappa_{i,j}(h) \equiv \binom{\tilde{\chi}(h)}{j-i} \bmod p^m
$$
for $h \in H$, where $\tilde{\chi} \colon H \to \Z$ is the lift of $\chi$ with image in $[0,p^n)$.  Allowing only these canonical maps in our defining systems, the ambiguity in the definition of the Massey product of $\chi_1, \ldots, \chi_q$ is reduced to the choice of the $\kappa_{i,q+1}$ with $2 \le i \le q-1$.  

We now pass to the setting of Galois cohomology so that we may interpret our
cohomology groups as groups of field elements (modulo powers) via Kummer
theory.  Suppose that $H$ is the Galois group of some Galois extension of fields $\Omega/F$, with $F$ of prime-to-$p$ characteristic containing the $p^n$th roots of unity $\mu_{p^n}$.  Then $\chi$ factors through the Galois group $Q$ of a Kummer extension $E/F$ of degree $p^n$.  The group $C_{q+1}$ is isomorphic to $(\Z/p^m\Z)[Q]/I_Q^q$ as a $(\Z/p^m\Z)[Q]$-module via conjugation by a lift to $T_{q+1}$ of the image of $Q$ in $T_{q+1}/C_{q+1}$.   We identify $\mu_{p^n}$ with $\Z/p^n\Z$ by a choice of primitive $p^n$th root of unity $\zeta_n$.  The Kummer dual 
$\Hom(C_{q+1},\mu_{p^n})$ of $C_{q+1}$ is then isomorphic to a submodule of $(\Z/p^m\Z)[Q]$.  Its generators as a group, labelled $D^{(i)}$ for nonnegative $i \le q-1$, satisfy
$$
	(\sigma-1)D^{(i)} = D^{(i-1)},
$$
with $D^{(0)}$ the norm element and $\sigma$ the generator of $Q$ which $\chi$ takes to $1$.  These ``derivative operators''
are studied in Section \ref{kummer}.

Note that $\chi$ (resp., $\chi_q$) corresponds to some $a \in F^{\times} - F^{\times p}$ (resp., $b \in F^{\times}$) that becomes a $p^n$th power (resp., $p^m$th power) in $\Omega$ by Kummer theory.  By the same principle, if we can realize $b$ as a norm of an element $y \in E^{\times}$, then the $D^{(i)}y$ with $i \le q-1$ can be used to define the maps $\kappa_{q-i,q+1}$, so long as $D^{(i)}y \in \Omega^{\times p^m}$ for $i \le q-2$.  We refer to a defining system so obtained as a {\em proper defining system} and the resulting $2$-cocycle as the ``$q$-fold Massey product of $a$ with $b$.''  
It is not hard to see that if $D^{(q-1)}y \in \Omega^{\times p^m}$, then $\rho$ is a homomorphism, so the class of the $2$-cocycle $\nu$ of \eqref{matrixcocycle} that defines the Massey product is trivial.  In general, $D^{(q-1)}y$ determines $\nu$ through a particular transgression map (Theorem \ref{masseytrans}).   Section \ref{embed} provides the tools needed to derive this, giving a cohomological discussion of embedding problems in a rather general context, comparing nonabelian cohomology and the sequence of base terms of the Hochschild-Serre spectral sequence. The indeterminacy in the $q$-fold Massey product of $a$ with $b$  is given by the group of $(q-1)$-fold Massey products of $a$ with $c \in F^{\times} \cap \Omega^{\times p^m}$.
 
\subsection{Massey products in number theory} 
 
We return to our earlier setting, beginning with a number field $K$ that
we now assume to contain the $p^n$th roots of unity.  Let $S$ be a set of primes of $K$ containing all primes above $p$ and any real places.  We consider Massey products for the maximal unramified outside $S$ extension $\Omega$ of $K$, using the construction of the previous subsection.  Theorem A arises by taking inverse limits  (increasing $n$) from a statement at the finite level, which in turn arises from a formula for the relevant Massey products.  This is the subject of Sections \ref{classgroups} and \ref{infiniteext}.

We quickly introduce the relevant notation.  Similarly to the previous subsection, let $a \in K^{\times} \cap \Omega^{\times p^n}$ with $a \notin
K^{\times p}$, let $b \in K^{\times} \cap \Omega^{\times p^m}$, let $L$ be the field given by adjoining a $p^n$th root $\alpha$ of $a$ to $K$, and let $\sigma$ be a generator of $L$ such that $\alpha^{\sigma-1} = \zeta_n$.  We use $P_{m,L/K}^{(k-1)}$ to denote the group of cohomology classes of all $k$-fold Massey products of $a$ with some $c \in K^{\times} \cap \Omega^{\times p^m}$ and set
$P_{m,L/K}^{(0)} = 0$.  The class of the $(k+1)$-fold Massey product of $a$ with $b$ is determined up to an element of $P_{m,L/K}^{(k-1)}$ when it admits a proper defining system.  We denote the class of the resulting Massey product modulo $P_{m,L/K}^{(k-1)}$ by $(a,b)^{(k)}_{p^m,K,S}$.

Kummer theory provides a natural injection
$$
  A_{K,S}/p^m A_{K,S} \hookrightarrow H^2(G_{K,S},\mu_{p^m}),
$$
where $A_{K,S}$ denotes the $p$-part of the $S$-class group of $K$.  The classes
of proper defining systems actually lie in the quotient of $A_{K,S}/p^m$ by the image of $P_{m,L/K}^{(k-1)}$ in this group (recalling our rigidification of $\mu_{p^n}$).  We prove a formula (Theorem \ref{pairformula}) for the Massey
products which essentially states the following.   If $(a,b)^{(k)}_{p^m,K,S}$ is defined, any proper defining system is given by the image
modulo $p^m$ of the ideal class of the norm of a fractional ideal
$\mathfrak{a}$ in the ring of $S$-integers of $K$ such that $\mathfrak{a}^{(\sigma-1)^k}$ is (nearly)
principal and has a generator of norm $b^{-1}$.  This formula is a generalization of
\cite[Theorem 2.4]{mcs}.
We use it to prove an analogue of our main result at the level of $K$ (Proposition \ref{isogroups}).

Passing up the cyclotomic tower, we can define Massey products at the level of $K_{\infty}$.  Fix a generator of the Tate module of $\mu_{p^{\infty}}$ so that we may identify it with $\zp$.  Write $K_n = K(\mu_{p^n})$ for each $n \ge 1$
(we no longer fix $n$ or assume $\mu_{p^n} \subset K$).  Assume, for simplicity of presentation, that $L_{\infty}$ is
defined by a sequence $a = (a_n)$ with
$$
  a_n \in K_n^{\times} \cap \Omega^{\times p^n},\ a_n \notin
  K_n^{\times p},\mr{\ and\ }
  a_{n+1}a_n^{-1} \in K_{n+1}^{\times p^n}.
$$
Let $r$ be such that $k < p^r(p-1)$.
Given a sequence $b$ 
of $p$-units $b_n \in K_n^{\times}$ for which the
images in $K_n/K_n^{\times p^{n-r}}$ are norm compatible
and the $(k+1)$-fold Massey products
$
  (a_n,b_n)_{p^{n-r},K_n,S}^{(k)}
$
admit proper defining systems, we may consider the inverse limit
of these products.  This defines an element $(a,b)_{K,S}^{(k)}$ 
of $\mc{A}_{K}$, the $(k+1)$-fold Massey product of $a$ with $b$.  We let $\mc{P}_{L/K}^{(k)}$ denote the submodule generated by all such inverse limits, varying $b$.  Theorem A (and more generally, Theorem \ref{mainthm}) then follows from its analogue at the finite level.

\begin{acknowledgments}
I thank Bill McCallum for our joint work on cup products
\cite{mcs} which served as the inspiration for this paper, Dick
Hain for introducing Massey products to me, Barry Mazur for his
encouragement and advice, and Mike Hopkins for an important
observation.  Section \ref{embed} is a condensed version of the
first chapter of my 1999 Ph.D. thesis \cite{me-thes}, and I thank
Spencer Bloch for advising me on this.  This research was
undertaken, in part, thanks to funding from an NSF Postdoctoral Research Fellowship, Harvard University, the Max Planck Institute for Mathematics,
and the Canada Research Chairs Program.
\end{acknowledgments}

\section{Galois embedding problems} \label{embed}

In this section, we study homomorphisms of the absolute Galois
group of a field to a finite group $\mc{H}$ and their lifts to a
larger finite group $\mc{G}$ with $\mc{H}$ as a quotient.  That
is, we study Galois embedding problems.  Our point of view differs
from those which we have been able to find in the literature
(see, for instance, \cite{ilf, koch, nsw}), so we hope that this
section will serve as a useful reference. The proofs are for the
most part straightforward and left to the reader.

We use the following notational conventions.  For a field $K$, we
let $G_K$ denote the Galois group of the separable closure of $K$
over itself.  For a Galois extension $L/K$, we let $G_{L/K} =
\Gal(L/K)$.

Now fix a field $K$.  We consider an exact sequence of finite
groups of the form
\begin{equation*}
  0 \to A \xrightarrow{\iota} \mc{G} \xrightarrow{\phi} \mc{H} \to 1,
\end{equation*}
where $A$ is abelian and view these groups as trivial
$G_K$-modules. We obtain a sequence of sets
\begin{equation} \label{homs}
  0 \to \Hom(G_K,A) \xrightarrow{\iota^.} \Hom(G_K,\mc{G})
  \xrightarrow{\phi^.} \Hom(G_K,\mc{H}),
\end{equation}
which is exact in the sense that those elements in the image of
one map are exactly those taken to the trivial homomorphism by the
next. Passing to nonabelian cohomology (e.g., \cite{serre}), we
have
\begin{equation} \label{nonab}
  0 \to H^1(G_K,A) \xrightarrow{\iota^*} H^1(G_K,\mc{G})
  \xrightarrow{\phi^*} H^1(G_K,\mc{H}).
\end{equation}
As $A$ is not necessarily central in $\mc{G}$, we cannot quite
extend the sequence \eq{nonab} to $H^2(G_K,A)$.  Instead, given a
homomorphism $\brho \in \Hom(G_K,\mc{H})$ we can twist $A$ by
$\brho$ so that the action of $\sigma \in G_K$ on $a \in A$ is now
given by
$$
   a^{\sigma} = f(\sigma) a f(\sigma)^{-1},
$$
where $f$ is any continuous function lifting $\brho$ to $\mc{G}$.
We denote this new module structure by $A_{\brho}$.  We then
obtain an element $\Delta(\brho) \in H^2(G_K,A_{\brho})$ by
lifting $\brho$ to some $f$ as above and defining the desired
$2$-cocycle by
$$
  a(\sigma_1,\sigma_2) = f(\sigma_1) f(\sigma_2) f(\sigma_1\sigma_2)^{-1}.
$$
The image of $a$ in cohomology is $\Delta(\brho)$, and this class
does not change for homomorphisms cohomologous to $\brho$.  In
particular, $\brho$ (resp., its class $[\brho]$) will be in the
image of $\phi^.$ (resp., $\phi^*$) if and only if $\Delta(\brho)$
is trivial.

\begin{remark}
  The class $\Delta(\brho)$ is the obstruction to lifting $\brho$ in a
  very real sense.  If $\Delta(\brho) = 0$, then $a$ is a coboundary,
  so we can choose $\kappa \colon G_K \to A$ with $d\kappa = -a$.
  Then, as one can easily check, $\kappa \cdot f$ is a homomorphism
  lifting $\brho$.
\end{remark}

We can give a description of $\Delta(\brho)$ in terms of group
extensions.  Consider the fiber product
\begin{equation*}
  \mc{G} \times_{\mc{H}} G_K = \langle (g,\sigma) \mid \phi(g) = \brho(\sigma)
  \co g \in \mc{G} \co \sigma \in G_K \rangle,
\end{equation*}
which is the pullback in the following diagram
\begin{equation*} \label{fp}
  \UseComputerModernTips \xymatrix{
  0 \ar[r] & A_{\brho} \ar[r] \ar@2{-}[d] & \mc{G} \times_{\mc{H}} G_K \ar[r]
  \ar[d] & G_K \ar[r] \ar[d]^{\brho} & 1 & (\dagger)\\
  0 \ar[r] & A \ar[r] & \mc{G} \ar[r]^{\phi} & \mc{H} \ar[r] & 1.
} \end{equation*}

\begin{lemma} \label{Delex}
  The class of $(\dagger)$ in $H^2(G_K,A_{\brho})$ is $\Delta(\brho)$.
\end{lemma}

Let $L$ denote the fixed field of the kernel of $\brho$.  Clearly,
we can also take the pullback
\begin{equation*}
  \UseComputerModernTips \xymatrix{
  0 \ar[r] & A_{\brho} \ar[r] \ar@2{-}[d] & \phi^{-1}(G_{L/K}) \ar[r]
  \ar[d] & G_{L/K} \ar[r] \ar@{^{(}->}[d]^{\brho} & 1
  & (\dagger\dagger)\\
  0 \ar[r] & A \ar[r] & \mc{G} \ar[r]^{\phi} & \mc{H} \ar[r] & 1.  }
\end{equation*}
The group extension $(\dagger\dagger)$ yields a class
$\delta(\brho) \in H^2(G_{L/K},A_{\brho})$.  However, this class
is not necessarily invariant under conjugation of $\brho$.  Hence
we do not define $\delta$ on the class $[\brho]$.  In particular,
we have by Lemma \ref{Delex} that
$$
  \Delta(\brho) = \Inf(\delta(\brho)),
$$
where $\Inf$ denotes the inflation from $G_{L/K}$ to $G_K$.

Next, we would like a suitable definition of the transgression map
in the Hochschild-Serre spectral sequence in terms of group
extensions.  This description can obviously be made in the general
setting of a profinite group $\mc{P}$, an open normal subgroup
$\mc{N}$, and a discrete $\mc{P}$-module $M$ upon which $\mc{N}$
acts trivially, but we restrict to the case of interest ($\mc{P} =
G_K$, and so forth) in this discussion.  So, consider the
transgression map
\begin{equation*}
  \Tra \colon H^1(G_L,A_{\brho})^{G_{L/K}} \to H^2(G_{L/K},A_{\brho}).
\end{equation*}
This map can be described as follows (e.g., \cite[Proposition
1.6.5]{nsw}). For a homomorphism $f \in
H^1(G_L,A_{\brho})^{G_{L/K}}$, extend $f$ to a continuous map $h
\colon G_K \to A_{\brho}$ satisfying
\begin{eqnarray} \label{transprop}
  h(\sigma\tau) = \sigma h(\tau) + h(\sigma) & \mr{and} &
  h(\tau\sigma) = h(\tau) + h(\sigma)
\end{eqnarray}
for all $\tau \in G_L$ and $\sigma \in G_K$. Then $dh$ induces
a well-defined 2-cocycle on $G_{L/K}$, the class of which is
$\Tra{f}$.

In terms of group extensions, we can express the transgression map
as follows. For $f \in H^1(G_L,A_{\brho})^{G_{L/K}}$, define $g:
G_L \to A_{\brho} \rtimes G_K$ by $g(\tau) = (f(\tau),\tau)$. Then
$g(G_L)$ is normal in $A_{\brho} \rtimes G_K$, so define $\mc{E} =
(A_{\brho} \rtimes G_K)/g(G_L)$. This yields the following
commutative diagram:
\begin{equation*}
  \UseComputerModernTips \xymatrix{
  0 \ar[r] & A_{\brho} \ar[r] \ar@2{-}[d] & A_{\brho} \rtimes G_K \ar[r]
  \ar@{>>}[d]
  & G_K \ar[r] \ar@{>>}[d] & 1 \\
  0 \ar[r] & A_{\brho} \ar[r] & \mc{E} \ar[r] & G_{L/K} \ar[r] & 1 & (*).  }
\end{equation*}
Alternatively, we can describe $\mc{E}$ as the pushout in the
following commutative diagram
$$
  \UseComputerModernTips \xymatrix{
  1 \ar[r] & G_L \ar[r] \ar[d]^{-f} & G_K \ar[r] \ar[d]
  & G_{L/K} \ar[r] \ar@{=}[d] & 1 \\
  0 \ar[r] & A_{\brho} \ar[r] & \mc{E} \ar[r] & G_{L/K} \ar[r] & 1 & (*).  }
$$

\begin{lemma} \label{tra}
  The equivalence class of $(*)$ as a group extension in $H^2(G_{L/K},A_{\brho})$
  is $\Tra{f}$.
\end{lemma}

We now link the two points of view discussed above: the study of
lifts of Galois representations and nonabelian cohomology, and the
study of group extensions and the transgression map.  Given $\rho
\in \Hom(G_K,\mc{G})$ with image $\brho \in \Hom(G_K,\mc{H})$ we
can define $\Lambda(\rho) \in H^1(G_L,A_{\brho})^{G_{L/K}}$ by
\begin{equation*}
  \Lambda(\rho)(\tau) = -\rho(\tau)
\end{equation*}
for $\tau \in G_L$.  Again, we do not define $\Lambda$ on the
class of $\rho$ unless $A$ is central in $\mc{G}$.

\begin{lemma} \label{tralam}
  We have $\Tra{\Lambda(\rho)} = \delta(\phi^.(\rho))$.
\end{lemma}

\begin{proof}
  Observe the following diagram:
  \begin{equation} \label{bigd}
    \UseComputerModernTips \xymatrix{
    1 \ar[r] & G_L \ar[r] \ar[d]^{-\Lambda(\rho)} & G_K \ar[r] \ar[d]
    & G_{L/K} \ar[r] \ar@{=}[d] & 1 \\
    0 \ar[r] & A_{\brho} \ar[r] \ar@{=}[d] & \mc{E} \ar[r] \ar@{.>}[d] &
    G_{L/K} \ar[r] \ar@{=}[d] & 1 \\
    0 \ar[r] & A_{\brho} \ar[r] \ar@2{-}[d] & \phi^{-1}(G_{L/K})
    \ar[r] \ar[d] & G_{L/K} \ar[r] \ar@{^{(}->}[d]^{\brho} & 1 \\
    0 \ar[r] & A \ar[r] & \mc{G} \ar[r]^{\phi} & \mc{H} \ar[r] & 1. }
  \end{equation}
  Note that $\phi^{-1}(G_{L/K}) = \mc{G} \times_{\mc{H}} G_{L/K}$.  We define
  $A_{\brho} \to \mc{G} \times_{\mc{H}} G_{L/K}$ by $a \mapsto (a,0)$ and $G_K \to
  \mc{G} \times_{\mc{H}} G_{L/K}$ by $\sigma \mapsto (\rho(\sigma),\bar{\sigma})$.
  These maps coincide on the images of $\tau \in G_L$ because
  $-\Lambda(\rho)(\tau) = \rho(\tau)$.  By the universal property of
  the pushout $\mc{E}$, we obtain a map $\mc{E} \to \phi^{-1}(G_{L/K})$
  for which \eq{bigd} commutes.  By the Five Lemma, it is an isomorphism.
\end{proof}

\begin{lemma} \label{timesp}
  Right multiplication by $\rho$ induces a well-defined map from
  the group of $1$-cocycles $Z^1(G_K,A_{\brho})$ to
  $\Hom(G_K,\mc{G})$ $($and, in fact, from
  $H^1(G_K,A_{\brho})$ to $H^1(G_K,\mc{G}))$.  Furthermore, for $k \in
  Z^1(G_K,A_{\brho})$ we have
  \begin{equation*}
    \Res{[k]} \, \cdot \, \Lambda(\rho) = \Lambda(k \, \cdot \, \rho),
  \end{equation*}
  where $\Res$ denotes restriction from $G_K$ to $G_L$.  Any lift of
  $\brho$ has the form $k \, \cdot \, \rho$ for some $k$.
\end{lemma}

\begin{lemma} \label{unique}
  Two lifts $\rho$ and $\rho'$ of $\brho$ satisfy $\Lambda(\rho) =
  \Lambda(\rho')$ if and only if $\rho' = t \rho$ for some $t \in
  Z^1(G_{L/K},A_{\brho})$.
\end{lemma}

Putting everything together, we have the following ``commutative
diagram'' in which certain of the maps, which we denote with
dashed arrows, must be taken only on certain elements (which map
to or come from $\brho$ or a lift $\rho$ of it) and in which
exactness holds only in the appropriate sense and where
appropriate.  We include it merely to summarize what we have
discussed above:
\begin{equation} \label{fake} \small
  \UseComputerModernTips \xymatrix@C=12pt{
  Z^1(G_{L/K},A_{\brho}) \ar@{^{(}->}[r] \ar@{>>}[d]
  & Z^1(G_K,A_{\brho}) \ar@{^{(}->}[r]^-{\iota^.} \ar@{>>}[d]
  & \Hom(G_K,\mc{G}) \ar[r]^-{\phi^.} \ar@{-->}[d]^{\Lambda}
  & \Hom(G_K,\mc{H}) \ar@{-->}[r]^-{\Delta} \ar@{-->}[d]^{\delta}
  & H^2(G_K,A_{\brho}) \ar@{=}[d] \\
  H^1(G_{L/K},A_{\brho}) \ar@{^{(}->}[r]^-{\Inf}
  & H^1(G_K,A_{\brho}) \ar[r]^-{\Res}
  & H^1(G_L,A_{\brho})^{G_{L/K}} \ar[r]^-{\Tra}
  & H^2(G_{L/K},A_{\brho}) \ar[r]^-{\Inf} & H^2(G_K,A_{\brho}). }
  \normalsize
\end{equation}

In particular, we obtain the following proposition which relates
liftings of Galois representations to group extensions.

\begin{proposition} \label{same}
  Let $\brho \colon G_K \to \mc{H}$, and let $L$ denote the fixed field of
  its kernel.  Then $\brho$ lifts to some $\rho \colon G_K \to \mc{G}$ if
  and only if $\delta(\brho)$ is in the image of the transgression
  map.  If $\delta(\brho) = \Tra{\xi}$, then $\rho$ may be chosen so
  that $\Lambda(\rho) = \xi$.  The choice of $\rho$ is unique up to
  left multiplication by a cocycle in $Z^1(G_{L/K},A_{\brho})$.
  Furthermore, $\rho$ will be surjective if and only if both $\brho$
  and $\xi$ are surjective.
\end{proposition}

\begin{remark}
  In the case that $A$ is central in $\mc{G}$, the module $A_{\brho}$ is just
  $A$ with trivial action, and all the maps we defined above pass to
  cohomology classes.  Hence we obtain another ``commutative
  diagram,'' similar to \eq{fake}, in which all the maps in the top
  row are defined on all elements:
  \begin{equation} \label{nice} \small
    \UseComputerModernTips \xymatrix@C=12pt{
    & 0 \ar[r] & H^1(G_K,A) \ar[r]^-{\iota^*} \ar@{=}[d] & H^1(G_K,\mc{G})
    \ar[r]^-{\phi^*} \ar@{-->}[d]^{\Lambda} & H^1(G_K,\mc{H}) \ar[r]^-{\Delta}
    \ar@{-->}[d]^{\delta} & H^2(G_K,A) \ar@{=}[d] \\
    0 \ar[r] & H^1(G_{L/K},A) \ar[r]^-{\Inf} & H^1(G_K,A) \ar[r]^-{\Res} &
    H^1(G_L,A)^{G_{L/K}} \ar[r]^-{\Tra} & H^2(G_{L/K},A) \ar[r]^-{\Inf} &
    H^2(G_K,A). }
  \end{equation}
\end{remark}

\begin{remark}
  The entire discussion already given in this section
  goes through with $G_K$ replaced by $G_{\Omega/K}$ for an
  algebraic Galois extension $\Omega/K$.
\end{remark}

In the case of $\Omega/K$, in view of the latter remark, we denote
the maps corresponding to $\Delta$ and $\Inf$ in \eq{fake} and
\eq{nice} with a subscript $\Omega$.  On the other hand, we use
$\Tra_{\Omega}$ to denote
$$
    \Tra_{\Omega} \colon H^1(G_{\Omega},A_{\brho})^{G_{\Omega/K}} \to
    H^2(G_{\Omega/K},A_{\brho})
$$
and $\Lambda_{\Omega}(\rho)$ to denote the restriction of
$\Lambda(\rho)$ to $G_{\Omega}$.

The following will be useful in Section \ref{massey}.

\begin{lemma} \label{omegalift}
  Let $\Omega$ be a Galois extension of $K$ containing $L$, and
  assume that $\brho$ lifts to a homomorphism $\rho \colon G_K \to \mc{G}$.
  The obstruction to lifting the map $\brho_{\Omega}$ induced by
  $\brho$ on $G_{\Omega/K}$ to an element of
  $\Hom(G_{\Omega/K},\mc{G})$ is given by the two equal quantities
  \begin{equation*} \label{equalobstr}
    \Delta_{\Omega}(\brho_{\Omega}) = \Tra_{\Omega}
    \Lambda_{\Omega}(\rho).
  \end{equation*}
\end{lemma}

\begin{proof}
  That $\Delta_{\Omega}(\brho_{\Omega})$ is the obstruction is immediate.
  By equation \eq{fake} for $G_{\Omega/K}$ and $G_K$, we have
  $$
    \Delta_{\Omega}(\brho_{\Omega}) = \Inf_{\Omega} \delta(\rho) =
    \Inf_{\Omega} \circ \Tra \Lambda(\rho) = \Tra_{\Omega} \Lambda_{\Omega}(\rho).
  $$
\end{proof}

\section{Kummer theory} \label{kummer}

We begin this section with an algebraic lemma, which we consider
in only slightly more generality than we shall need.  Let $G$ be a
finite group, and let $R$ be $\Z/N\Z$ for some $N$. Let $\iota
\colon R[G] \to R[G]$ denote the linear map determined by the
inversion map $g \mapsto g^{-1}$ on $G$. Given a left ideal $J$ of
$R[G]$, we may consider the right ideal $J^{\iota}$ of $R[G]$. Let
$J^{\perp}$ denote the left 
annihilator of $J^{\iota}$, which is a left
ideal of $R[G]$.

\begin{lemma} \label{duality}
Let $J$ be a left ideal in $R[G]$.  Then
$$
  \phi \colon \Hom(R[G]/J, \Q/\Z) \to J^{\perp},
$$
with $\phi(f) = \sum_{g \in G} f(g)g$, is an isomorphism of left
$R[G]$-modules.
\end{lemma}

\begin{proof}
  First, we remark that
  $$
    \phi(hf) = \sum_{g \in G} f(h^{-1}g)g = h \sum_{g \in G}
    f(g)g.
  $$
  Hence, $\phi$ is a homomorphism of left $R[G]$-modules as a map into
  $R[G]$.  Since $Q = \Hom(R[G]/J, \Q/\Z)$ is annihilated by
  $J^{\iota}$, the image of $\phi$ is contained in $J^{\perp}$.
  Furthermore, we claim that there exists a homomorphism of left
  $R[G]$-modules,
  $\psi \colon J^{\perp} \to Q$, given by
  $$
    a = \sum_{g \in G} a_g g \mapsto (g \mapsto a_g).
  $$
  To see this, we remark that it exists as a map to
  $\Hom(R[G],\Q/\Z)$, and for $r = \sum_{g \in G} r_g g \in J$, 
  we have
  $$
    \psi(a)(r) = \sum_{g \in G} r_g \psi(a)(g) = \sum_{g \in G}
    r_g a_g,
  $$
  which is the coefficient of $1$ in $\iota(r)a = 0$.  Then $\psi$
  is clearly the inverse map to $\phi$.
\end{proof}

\begin{remark}
If $J'$ is a left ideal containing $J$, then we have a commutative
diagram
$$
  \xymatrix{
  0 \ar[r] & \Hom(R[G]/J',\Q/\Z) \ar[r] \ar[d] & \Hom(R[G]/J,\Q/\Z) \ar[r] \ar[d] &
  \Hom(J'/J,\Q/\Z) \ar[r] \ar[d] & 0 \\
  0 \ar[r] & (J')^{\perp} \ar[r] & J^{\perp} \ar[r] &
  J^{\perp}/(J')^{\perp} \ar[r] & 0}
$$
in which the vertical arrows are isomorphisms.  We also have
$(J^{\perp})^{\perp} = J$.
\end{remark}

We will be applying Lemma \ref{duality} in the setting of Kummer
theory. Let $K$ be a field of characteristic not equal $p$
containing the group $\mu_{p^n}$ of $p^n$th roots of unity. Let
$L/K$ be a Galois extension, denote its Galois group by $G$, and
consider a subgroup $\Phi$ of $L^{\times}$ which is a left
$\Z[G]$-module.  Then the field $M$ obtained by adjoining all
$p^n$th roots of elements of $\Phi$ to $L$ is Galois over $K$. Let
$$
  \bar{\Phi} = \Phi/(\Phi \cap L^{\times p^n}).
$$
Kummer theory provides a (noncanonical) isomorphism
$$
  G_{M/L} \cong \Hom(\bar{\Phi},\Q/\Z).
$$
If $\bar{\Phi}$ is generated by one element as a module, we may
apply Lemma \ref{duality} with $R = \Z/p^n\Z$ to obtain a second
description of this Galois group.

Let us consider the case that $G$ is cyclic of order $p^n$ with
generator $\sigma$.  Let
\begin{equation*} \label{Dkdef}
  D_{\sigma}^{(k)} =
  (-1)^k\sum_{i=0}^{p^n-1}\binom{i}{k}\sigma^{i-k}
  \in \Z[G]
\end{equation*}
for $0 \le k \le p^n-1$.  When no confusion about the choice of
generator can result, we denote $D_{\sigma}^{(k)}$ simply by
$D^{(k)}$.  Letting $\Z[G]$ act on $L^{\times}$ in the natural
way, $D^{(0)}$ acts as $\Norm$, where $\Norm$ denotes the norm for
the extension $L/K$.  We leave the following lemma to the reader
to verify.

\begin{lemma} \label{Dk}
  Let $k$ be a positive integer with $k \le p^n-1$.  Then $D^{(k)}$
  satisfies
  $$
    (\sigma-1)D^{(k)} = D^{(k-1)} + (-1)^k \binom{p^n}{k}\sigma^{-k}
    \equiv D^{(k-1)} \bmod p^{n-v_p(k)}\Z[G].
  $$
  Therefore, letting $r_k = \max \{r \mid k \ge p^r \}$, we have
  $$
    (\sigma-1)^j D^{(k)} = D^{(k-j)} + (-\sigma)^{j-1-k} \sum_{i=0}^{j-1}
    \binom{p^n}{k-i} (\sigma^{-1}-1)^{j-1-i} \equiv D^{(k-j)} \bmod
    p^{n-r_k}\Z[G]
  $$
  for $1 \le j \le k$.
\end{lemma}

Let $I_G$ denote the augmentation ideal of $R[G]$ for a given
group $G$ and coefficient ring $R$.

\begin{lemma} \label{auggen}
  Let $k$ and $m$ be nonnegative integers with $m + r_k \le n$
  $($that is, $k \le p^{n-m+1}-1)$, and set $R = \Z/p^m\Z$.
  Then the ideal $(I_G^{k+1})^{\perp}$ of $R[G]$
  is generated by $D^{(k)}$.
\end{lemma}

\begin{proof}
  Note that $(I_G^t)^{\iota} = I_G^t$ for all $t$.
  For $k = 0$, the statement of the lemma is obvious.
  We have the inclusion $(I_G^k)^{\perp} \subset (I_G^{k+1})^{\perp}$.
  Since $I_G^k/I_G^{k+1}$ is cyclic of order dividing $p^m$, so is
  $(I_G^{k+1})^{\perp}/(I_G^k)^{\perp}$.  It is then clearly
  generated by the element $D^{(k)}$ of order $p^m$, as
  $(\sigma-1)D^{(k)} = D^{(k-1)}$ in $R[G]$.
\end{proof}

It is not hard to pass from Lemma \ref{auggen} to the general case
of any $k$, $m$, and $n$, but we shall not require this.  On the
other hand, in Section \ref{classgroups}, we shall require the
following lemma that describes which $D^{(k)}$ become trivial in
the group rings of quotients of $G$.

\begin{lemma} \label{trivimage}
  Let $k$ be a nonnegative integer and
  $m$ and $n$ be positive integers with $k \le p^{n-m+1}-1$.
  Let $G$ be a cyclic group of order $p^n$ with generator $\sigma$,
  and let $H$ be its quotient of order $p^{n-m}$.  Let
  $\phi \colon R[G] \to R[H]$ denote the
  natural map on group rings, for $R = \Z/p^m\Z$.  Then we have
  $$
    \phi(D_{\sigma}^{(k)}) =
    \begin{cases}
      0 & \mr{if\ } k < p^{n-m}(p-1),\\
      p^{m-1} D_{\phi(\sigma)}^{(k-p^{n-m}(p-1))} & \mr{if\ } k \ge p^{n-m}(p-1).
    \end{cases}
  $$
\end{lemma}

\begin{proof}
  The proof we give is combinatorial.
  We have
  \begin{equation} \label{phidk}
    \phi(D_{\sigma}^{(k)}) = (-1)^k\sum_{i=0}^{p^{n-m}-1}\sum_{j=0}^{p^m-1}
    \binom{i+jp^{n-m}}{k}\phi(\sigma)^{i-k}.
  \end{equation}
  Note that
  \begin{equation} \label{binomsum}
    \binom{i+jp^{n-m}}{k} =
    \sum_{a=0}^{k} c_{i,k}(a) (jp^{n-m})^a,
  \end{equation}
  where $c_{i,k}(0) = \binom{i}{k}$ and, for $a > 0$,
  \begin{equation} \label{coeff}
    c_{i,k}(a) =  \sum_{0 \le b_1 < \ldots < b_a < k} \frac{1}{k!}
    \prod_{\substack{d = 0 \\ d \notin \{ b_1, \ldots, b_a
    \}}}^{k-1} (i-d).
  \end{equation}
  Using the facts that $i < p^{n-m}$ and $k < p^{n-m+1}$,
  one can obtain a lower bound on the $p$-adic valuations of
  the terms in the sum in \eq{coeff} which yields
  \begin{equation} \label{valuation}
    p^{a(n-m)}c_{i,k}(a) \in
    \begin{cases}
      \Z_p & \mr{if\ } a \le p-1, \\
      p\Z_p & \mr{if\ } a \ge p.
    \end{cases}
  \end{equation}
  Furthermore, one may check that
  \begin{equation} \label{powersum}
    \sum_{j=0}^{p^m-1} j^a \equiv
    \begin{cases}
      -p^{m-1} \bmod p^m & \mr{if\ } a \ge 1 \mr{\ and\ } a \equiv 0 \bmod p-1,\\
      0 \bmod p^m & \mr{if\ } a = 0 \mr{\ or\ } a \not\equiv 0 \bmod p-1
    \end{cases}
  \end{equation}
  unless $p = 2$, $m \ge 2$, $a$ is odd and $a \ge 3$, in which case
  the sum is $0 \bmod p^m$.
  Note that equations \eq{binomsum}, \eq{valuation}, and \eq{powersum}
  reduce the computation to that of $c_{i,k}(p-1)$.

  Now let us consider the specific case that $k = p^{n-m+1}-1$.
  One checks in a manner similar to that of \eq{valuation} that
  \begin{equation} \label{specific}
    p^{(p-1)(n-m)}c_{i,p^{n-m+1}-1}(p-1) \equiv
    \begin{cases}
    0 \bmod p & \mr{if\ } i \le p^{n-m}-2, \\
    -1 \bmod p & \mr{if\ } i = p^{n-m}-1.
    \end{cases}
  \end{equation}
  Putting equations \eq{binomsum}, \eq{powersum}, and
  \eq{specific}
  together, we obtain
  $$
    \sum_{j=0}^{p^m-1} \binom{i+jp^{n-m}}{p^{n-m+1}-1}
    \equiv \begin{cases}
      0 \bmod p^m & \mr{if\ } i \le p^{n-m}-2, \\
      p^{m-1} \bmod p^m & \mr{if\ } i = p^{n-m}-1.
    \end{cases}
  $$
  The result for $k = p^{n-m+1}-1$ now follows by recalling equation
  \eq{phidk}, since
  $$
    p^{m-1}D_{\phi(\sigma)}^{(p^{n-m}-1)} \equiv p^{m-1} \bmod p^m\Z[H]
  $$
  The result for all $k$ then follows by applying
  Lemma \ref{Dk}.
\end{proof}

It would be nice to have a more intuitive proof of Lemma \ref{trivimage}.
In any case, we can derive the following from it.

\begin{lemma} \label{projform}
  Let $k$, $s$, $m$, and $n$ be nonnegative integers satisfying $1 \le m \le n$,
  $s \le n$, and $k \le p^{n-m}(p-1)$.  Let $G$ be a cyclic group of order
  $p^n$ with generator $\sigma$, and let $H$ be its quotient of order
  $p^{n-s}$.  Let $R = \Z/p^m\Z$, and let $\phi \colon R[G] \to R[H]$ be the
  natural map on group rings.
  Let $J$ be the ideal of $R[H]$ generated by $\phi(\sigma)^{p^{n-m}}-1$.
  Then
  $$
    \phi((I_G^k)^{\perp})^{\perp} \subseteq I_H^k + J^{\perp}.
  $$
\end{lemma}

\begin{proof}
  The cases in which $s = 0$ and $s \ge m$ are trivial, so we assume
  that $1 \le s < m \le n$.  If $k = 0$, then again
  the statement is trivial.  Assume that we know the result for $k-1$ with
  $k \le p^{n-m}(p-1)$.  Then we have that
  \begin{equation*} \label{phicont}
    \phi((I_G^k)^{\perp})^{\perp} \subseteq
    \phi((I_G^{k-1})^{\perp})^{\perp} \subseteq I_H^{k-1} + J^{\perp}.
  \end{equation*}
  Take $x \in \phi((I_G^k)^{\perp})^{\perp}$, so that
  \begin{equation} \label{cong}
    x \equiv a (\phi(\sigma)-1)^{k-1}  \bmod (I_H^k + J^{\perp}).
  \end{equation}
  for some $a \in R$.  Now $x \phi(D_{\sigma}^{(k-1)}) = 0$ by
  Lemma \ref{auggen}, so noting Lemma \ref{Dk}, we have
  \begin{equation} \label{ad}
    a \phi(D_{\sigma}^{(0)}) = ap^s D_{\phi(\sigma)}^{(0)}
    \in J^{\perp}\phi((I_G^k)^{\perp}).
  \end{equation}
  We define a generator $\bar{N}$ of $J^{\perp}$ by
  $$
    \bar{N} = \sum_{i=0}^{p^{m-s}-1} \phi(\sigma)^{p^{n-m}i}.
  $$
  Consider the quotient $H'$ of $G$ of order $p^{n-m}$, and
  let $\phi'$ denote the map on group rings.
  Since
  $$
    D_{\phi(\sigma)}^{(0)} = \bar{N} \sum_{i=0}^{p^{n-m}-1} \phi(\sigma)^i,
  $$
  the condition of \eqref{ad} may be rewritten as
  \begin{equation*} \label{apd}
    ap^s D_{\phi'(\sigma)}^{(0)} \in \phi'((I_G^k)^{\perp}).
  \end{equation*}
  Now, by Lemma \ref{trivimage}, we have $\phi'((I_G^k)^{\perp}) = 0$,
  since $k \le p^{n-m}(p-1)$.
  Hence, we must have $a = p^{m-s}b$ for some $b \in R$.
  Going back
  to equation \eqref{cong}, we see that
  $$
    x \equiv bp^{m-s}(\phi(\sigma)-1)^{k-1} \equiv
    b\bar{N}(\phi(\sigma)-1)^{k-1}
    \equiv 0 \bmod (I_H^k+J^{\perp}),
  $$
  finishing the proof.
\end{proof}

\section{Massey products} \label{massey}

We now fix a field $K$ of characteristic not equal to $p$
containing $\mu_{p^n}$ for some $n \ge 1$.  We also fix a Galois
extension $\Omega/K$.  As before, for $k \ge 1$, let $r_k$ be the
largest integer such that $p^{r_k} \le k$.  Fix $k \le p^n-1$ and
$m \le n - r_k$, and consider an array of the form
$$
  C = ( C_{i,j} = \mu_{p^m}^{\otimes j-i} \mid 1 \le i < j \le
  k+2 ).
$$
This forms a multiplicative system as in \cite[p.\
183]{dwyer}, in that there are multiplication maps
$$
  C_{i,j} \otimes C_{j,k} \to C_{i,k}.
$$

Let $a \in K^{\times} \cap \Omega^{\times p^n}$ and $b \in
K^{\times} \cap \Omega^{\times p^m}$.  We will be interested in
the $(k+1)$-fold Massey product of $k$-copies of $a$ with $b$,
denoted $(a,b)_{p^m,\Omega/K}^{(k)}$ when defined, which is to say
it has a defining system described as follows.  A defining system
is a set of continuous maps
$$
  \kappa_{i,j} \colon G_{\Omega/K} \to C_{i,j}
$$
for $1 \le i < j \le k+2$ with $(i,j) \neq (1,k+2)$, that satisfy
\begin{enumerate}
\item[a.] $\kappa_{i,i+1}$ is the Kummer character associated to a $p^m$th
root of $a$ for all positive $i \le k$,
\item[b.] $\kappa_{k+1,k+2}$ is the Kummer character associated to a $p^m$th
root of $b$, and
\item[c.] for any $i$ and $j$ as above and 
$\sigma_1, \sigma_2 \in G_{\Omega/K}$, we have
$$
  \kappa_{i,j}(\sigma_1\sigma_2) = \kappa_{i,j}(\sigma_1) +
  \kappa_{i,j}(\sigma_2) +
  \sum_{l=1}^{j-i-1} \kappa_{i,l+i}(\sigma_1)\kappa_{l+i,j}(\sigma_2).
$$
\end{enumerate}
The $(k+1)$-fold Massey product $(a,b)^{(k)}_{p^m,\Omega/K}$ with
respect to such a defining system is the class in
$H^2(G_{\Omega/K},\mu_{p^m}^{\otimes k+1})$ of the $2$-cocycle
\begin{equation} \label{Masseyco}
  (\sigma_1,\sigma_2) \mapsto
  \sum_{i=2}^{k+1} \kappa_{1,i}(\sigma_1)\kappa_{i,k+2}(\sigma_2).
\end{equation}
Therefore, when the Massey product is defined, it is defined only
up to a choice of defining system (in fact, modulo some subgroup
of $H^2(G_{\Omega/K},\mu_{p^m}^{\otimes k+1})$). Later, we shall
restrict to a subset of all defining systems in order to reduce
this ambiguity (and, consequently, the set of Massey products
which are considered to be defined).

We may interpret defining systems for Massey products in terms of
matrices. Throughout, we fix a choice of primitive $p^n$th root of
unity $\zeta$, which we use to make an identification $\mu_{p^n}
\cong \Z/p^n\Z$.  We shall later show the resulting Massey product
to be independent of this choice. We envision the multiplicative
system as defining the group of upper triangular unipotent
matrices $\mc{T}_m^{(k)}$ in $GL_{k+2}(\Z/p^m\Z)$.  Let
$\mc{Z}_m^{(k)}$ denote the center of this group (the cyclic
subgroup generated by the matrix $Z = I + E_{1,k+2}$, where
$E_{i,j}$ is the matrix which is zero in all but the $(i,j)$th
entry, which is $1$). A defining system for
$(a,b)^{(k)}_{p^m,\Omega/K}$ then provides a homomorphism
$$
  \brho \colon G_{\Omega/K} \to \mc{T}_m^{(k)}/\mc{Z}_m^{(k)},
$$
as in \eq{masseymatrix}, taking $\chi_i = \kappa_{i,i+1}$.  By
definition of $\Delta_{\Omega}$ as in Section \ref{embed} with $A
= \mc{Z}^{(k)} \cong \Z/p^m\Z$, the cocycle defining
$\Delta_{\Omega}(\brho)$ is, under the appropriate identification
(taking $Z$ to $1$), the negative of the cocycle in \eq{Masseyco}
for the Massey product relative to the given defining system.

Let $\mc{G}_m^{(k)}$ be the subgroup of $\mc{T}_m^{(k)}$ generated
by the two elements $X = I + \sum_{i=1}^k E_{i,i+1}$ and $Y = I +
E_{k+1,k+2}$.  The quotient $H_m$ of $\mc{G}_m^{(k)}$ by the
normal subgroup $\mc{N}_m^{(k)}$ generated by $Y$ is cyclic of
order $p^{m+r_k}$, generated by the image of $X$.  We remark that
$\mc{G}_m^{(k)} \cong \mc{N}_m^{(k)} \rtimes H_m$.  Letting $G =
H_{n-r_k}$, which has order $p^n$, we may view $H_m$ canonically
as a quotient of $G$ for any fixed $m$ with $m + r_k \le n$.

For $i$ with $0 \le i \le k$, we let $Y_i = I + E_{k+1-i,k+2}$,
so that
$$
  Y_i = [X,[X,\ldots[X,Y]\ldots ]],
$$
where the commutator is of length $i+1$ (using the convention $[X,Y]
= XYX^{-1}Y^{-1}$). The group $\mc{N}^{(k)}$ is isomorphic to
$(\Z/p^m\Z)[G]/I_G^{k+1}$ as a $(\Z/p^m\Z)[G]$-module, under the
isomorphism taking $Y$ to $1$.

Let $\alpha^{p^n} = a$, and set $L = K(\alpha)$.  We assume for
simplicity that $[L:K] = p^n$ throughout.  Let $\sigma$ generate
$G_{L/K}$, and let $\zeta$ be such that $\sigma \alpha = \zeta
\alpha$.  For any $t \le n$, let $\zeta_t = \zeta^{p^{n-t}}$.
Assume that the norm residue symbol $(a,b)_{p^n,K}$ is trivial,
which is to say that there exists $y \in L^{\times}$ with $N_{L/K}
y = b$.  We identify $G_{L/K}$ with $G$ under the map taking
$\sigma$ to the image of $X \in \mc{G}_{n-r_k}^{(k)}$.  Let $\Phi$
denote the $\Z[G]$-submodule of $L^{\times}$ generated by $y$. Let
$M$ denote the Kummer extension of $L$ given by adjoining the
$p^m$th roots of all elements of $\Phi$. In the notation of
Section \ref{kummer}, we have an identification $\bar{\Phi} \cong
(\Z/p^m\Z)[G]/J$ of $(\Z/p^m\Z)[G]$-modules taking the image of
$y$ to the image of $1$, for some ideal $J$.  Lemma \ref{duality}
then implies that $G_{M/L} \cong J^{\perp}$ (with $R = \Z/p^m\Z$,
under the isomorphism defined by $\zeta_m$). For $j \le k$, let
$M_j$ denote the subextension corresponding to the submodule
$((I_G^{j+1})^{\perp}+J)/J$ of $\bar{\Phi}$.  Also, let $c_j =
D^{(j)}y$ and $\eta_j^{p^m} = c_j$. By Lemma \ref{auggen}, we have
$M_j = L(\eta_0,\ldots,\eta_j)$. We remark that $G_{M_j/L}$ has
the form $J^{\perp}/(I_G^{j+1} \cap J^{\perp})$ as a
$(\Z/p^m\Z)[G]$-module.

\begin{lemma} \label{defining}
  For positive integers $k$ and $m$ with $m + r_k \le n$,
  there is a homomorphism
  $$ \rho^{(k)} \colon G_{M_k/K} \to \mc{G}_m^{(k)} $$
  under which
  $G_{M_k/L}$ is considered as a subgroup of $\mc{N}_m^{(k)}$
  via the injection
  \begin{equation} \label{submod}
    J^{\perp}/(I_G^{k+1} \cap J^{\perp}) \to (\Z/p^m\Z)[G]/I_G^{k+1}
  \end{equation}
  and for which, fixing a lift $\tilde{\sigma}$
  of $\sigma$ to $G_{M_k/K}$, we have
  $$
    \rho^{(k)}(\tilde{\sigma})
    = X + tE_{k+1,k+2},
  $$
  where $\zeta_m^t = \lambda(\tilde{\sigma})$, where $\lambda$
  is the Kummer character associated to a $p^m$th root of $b$.
\end{lemma}

\begin{proof}
  We must check that $\rho^{(k)}$, defined on $G_{M_k/L}$ and
  $\tilde{\sigma}$ as in the lemma, extends to a homomorphism
  on $G_{M_k/K}$ to the semidirect product $\mc{G}_m^{(k)}$.
  To see this, we begin by noting that
  \begin{equation} \label{action}
    Y^{(\rho^{(k)}(\tilde{\sigma})-1)^j} = Y_j
  \end{equation}
  for each $j$ with $0 \le j \le k$.

  We consider an element $e$ of $G_{M_k/L}$ as
  an element of $(\Z/p^m\Z)[G]/I_G^{k+1}$
  via the injection in \eq{submod}.
  Note that the norm element $N_G \in (\Z/p^m\Z)[G]$ satisfies
  $$
    N_G \equiv p^{m-1}(\sigma-1)^{p^{n-m+1}-1} \bmod
    I_G^{k+1}.
  $$
  Thus, unless $e = 0$, or $k = p^{n-m+1}-1$ and
  $e$ is a multiple of $p^{m-1}(\sigma-1)^{p^{n-m+1}-1}$,
  we have that $\sigma \tilde{e} \neq \tilde{e}$
  for any lift of $e$ to
  $$
    G_{M/L} \cong J^{\perp} \le (\Z/p^m\Z)[G].
  $$
  Therefore, if $k \le p^{n-m+1}-2$, then $\tilde{\sigma}^{p^n} = 1$ and
  $$
    G_{M_k/K} \cong G_{M_k/L} \rtimes G,
  $$
  so \eq{action} suffices to prove the lemma in this case.

  Now take $k = p^{n-m+1}-1$, and set $\rho = \rho^{(k)}$.  We are
  reduced to checking that
  $\rho(\tilde{\sigma}^{p^n}) = \rho(\tilde{\sigma})^{p^n}$.
  As $\tilde{\sigma}^{p^n}$ is an element of $G_{M/L}$ fixed by
  the action of $G_{L/K}$, we have seen that
  $\rho(\tilde{\sigma}^{p^n})$ is a power of
  $Y_{p^{n-m+1}-1}$.
  Note also that
  $$
    \rho(\tilde{\sigma})^{p^n} = Y_{p^{n-m+1}-1}^{p^{m-1}t},
  $$
  so we need only check that
  $$
    \lambda(\tilde{\sigma})^{p^{m-1}} =
    \lambda'(\tilde{\sigma}^{p^n}),
  $$
  where $\lambda'$
  denotes the Kummer
  character associated to a $p^m$th root of $c_{p^{n-m+1}-1}$.
  This follows from the following two facts in the appropriate group
  rings:
  $$
    \tilde{\sigma}^{p^n}-1 = (\tilde{\sigma}-1)\sum_{i=0}^{p^n-1}
    \tilde{\sigma}^i
  $$
  and
  $$
    N_G D^{(p^{n-m+1}-1)} = p^{m-1} N_G.
  $$
\end{proof}

Lemma \ref{defining} has the following consequence for defining
systems of Massey products using their definition in terms of
matrices. For $x \in \Omega^{\times}$, let $[x]$ denote the
cohomology class in $H^1(G_{\Omega},\mu_{p^m})$ associated to $x$
by Kummer theory.

\begin{proposition} \label{definingcor}
  Let $\rho^{(k)}$ be the homomorphism defined in Lemma
  \ref{defining}.
  If $c_{k-1} \in \Omega^{\times p^m}$, then
  the induced homomorphism
  $$
    \brho_{\Omega}^{(k)} \colon G_{\Omega/K} \to
    \mc{G}_m^{(k)}/\mc{Z}_m^{(k)}
  $$
  provides a defining system for the $(k+1)$-fold Massey
  product $(a,b)^{(k)}_{p^n,\Omega/K}$ with class
  \begin{equation*} \label{cocycles}
    \Tra_{\Omega} [c_k] \otimes \zeta_m^{\otimes k}
    \in H^2(G_{\Omega/K},\mu_{p^m}^{\otimes k+1}).
  \end{equation*}
  Given $y$ and $\zeta$, this class does not depend upon
  the choice of lift $\tilde{\sigma}$ of $\sigma$.
\end{proposition}

\begin{proof}
  The first statement follows by definition, using the
  identification given by the choice of $\zeta$, since the
  assumption on $c_{k-1}$ implies that $M_{k-1} \subseteq \Omega$.
  The class of the defining system is
  $-\Delta_{\Omega}(\brho^{(k)}) \otimes \zeta_m^{\otimes k}$.
  By Lemma \ref{omegalift} (with $L = M_{k-1}$), we have
  $$
    \Delta_{\Omega}(\brho_{\Omega}^{(k)}) =
    \Tra_{\Omega} \Lambda_{\Omega}(\rho),
  $$
  where $\rho$ denotes the inflation of $\rho^{(k)}$ to $G_K$.
  By definition, $\Lambda_{\Omega}(\rho)$
  is the negative of the class of $\rho|_{G_{\Omega}}$.
  By the injection in \eq{submod}, this class is equal to $[c_k]$.
  That this class does not depend upon the choice of lift
  $\tilde{\sigma}$ is clear.
\end{proof}

We make the following definitions for $m \le n-r_k$. If there
exists $y \in L^{\times}$ with $N_{L/K} y = b$ and $c_{k-1} \in
\Omega^{\times p^m}$, then the defining system provided by $y$ as
in Proposition \ref{definingcor} is called a proper defining
system.  Let $U_{m,L/K}^{(k)}$ denote the group of $b \in K^{\times}
\cap \Omega^{\times p^m}$ for which the Massey product
$(a,b)_{p^m,\Omega/K}^{(k)}$ has a proper defining system. Let
$P_{m,L/K}^{(0)}$ be the trivial group, and for $k \ge 1$, let
$P_{m,L/K}^{(k)}$ denote the group of classes in
$H^2(G_{\Omega/K},\mu_{p^m})$ which are twists by
$\zeta_m^{\otimes -k}$ of proper defining systems of
$(a,b)^{(k)}_{p^m,\Omega/K}$ for any $b \in U_{m,L/K}^{(k)}$. By
Proposition \ref{definingcor}, we may write this more precisely as
$$
  P_{m,L/K}^{(k)} = \langle \Tra_{\Omega}\,[D^{(k)}y] \mid
  y \in L^{\times}, D^{(k-1)}y \in \Omega^{\times p^m} \rangle.
$$
When we wish to indicate the
dependence of these definitions on $\Omega$, we will write
$U_{m,\Omega/L/K}^{(k)}$
and $P_{m,\Omega/L/K}^{(k)}$ in place of $U_{m,L/K}^{(k)}$ and
$P_{m,L/K}^{(k)}$, respectively.

Restricting to the set of proper defining systems (for a fixed
$n$) reduces the ambiguity in the class of the Massey product
$(a,b)^{(k)}_{p^m,\Omega/K}$ to the choice of an $x \in
L^{\times}$ with $D^{(k-2)}x \in \Omega^{\times p^m}$, as we may
replace $y$ with $N_{L/K}y = b$ by $yx^{\sigma-1}$. This changes
the resulting class by a twist of the class associated with the
proper defining system of $(a,N_{L/K}x)_{p^m,\Omega/K}^{(k-1)}.$
In particular, $P_{m,L/K}^{(k)}$ contains $P_{m,L/K}^{(k-1)}$.

From now on, we take the Massey product
$(a,b)^{(k)}_{p^m,\Omega/K}$ for $b \in U_{m,L/K}^{(k)}$ to be the
unique element of
$$
  H^2(G_{\Omega/K},\mu_{p^m}^{\otimes k+1})/(P_{m,L/K}^{(k-1)}
  \otimes \mu_{p^m}^{\otimes k})
$$
associated to any proper defining system.  (We will show this to
be independent of the choice of $\zeta$ in Proposition
\ref{welldef}.  We will also examine its dependence upon $n$,
which we exclude from the notation.) With these definitions, the
following theorem, which is the main result of this section, is a
direct corollary of Proposition \ref{definingcor}.

\begin{theorem} \label{masseytrans}
  Let $a \in K^{\times} \cap \Omega^{\times p^n}$, let $L =
  K(\alpha)$ for $\alpha \in \Omega^{\times}$ with $\alpha^{p^n} =
  a$, and assume that $[L:K] = p^n$.
  Let $\sigma$ denote a
  generator of $G_{L/K}$ with $\sigma\alpha = \zeta\alpha$.
  Let $k$ and $m$ be positive integers with $m+r_k \le n$.
  Let $b \in U_{m,L/K}^{(k)}$, and let $y \in L^{\times}$ be such that
  $N_{L/K} y = b$ and $D_{\sigma}^{(k-1)}y \in \Omega^{\times p^m}$.
  Then the $(k+1)$-fold Massey product satisfies
  $$
      (a,b)^{(k)}_{p^m,\Omega/K} \otimes \zeta_m^{\otimes -k}  =
      \Tra_{\Omega} [D_{\sigma}^{(k)} y] \pmod{P_{m,L/K}^{(k-1)}},
  $$
  where $\zeta_m = \zeta^{p^{n-m}}$.
\end{theorem}

We now verify that our definition of the Massey product is
independent of the original choice of $\zeta$.

\begin{proposition} \label{welldef}
  The Massey product $(a,b)^{(k)}_{p^m,\Omega/K}$
  is independent of the choice of $\zeta$.  Furthermore, $P_{m,L/K}^{(k)}$
  is independent of the choice of the element $a \in K$
  defining $L$.
\end{proposition}

\begin{proof}
  If we replace $\zeta$ by $\zeta^j$ for some $j$ prime to $p$, then
  by its definition, $\sigma$ is replaced by $\sigma^j$ as well.
  By Lemma \ref{auggen}, both $D^{(k)}$ and $D_j^{(k)} = D_{\sigma^j}^{(k)}$
  have image generating
  $(I_G^{k+1})^{\perp}/(I_G^k)^{\perp}$ for $R = \Z/p^m\Z$, and
  so they agree modulo
  $(I_G^k)^{\perp}$ up to a constant.  To compute this constant,
  note $(\sigma-1)^k D^{(k)} = N_G$ and
  $$
    N_G = (\sigma^j-1)^k D_j^{(k)} = \left( \sum_{i=0}^{j-1} \sigma^i
    \right)^k \cdot (\sigma-1)^k D_j^{(k)} = j^k (\sigma-1)^k
    D_j^{(k)}.
  $$
  Since $D^{(k-1)}$ generates $(I_G^k)^{\perp}$ as a $(\Z/p^m\Z)[G]$-module
  by Lemma \ref{auggen}, we then have that
  $$
    \Tra_{\Omega} [c_k] \equiv j^k \Tra_{\Omega} [D_j^{(k)} y] \bmod
    P_{m,L/K}^{(k-1)},
  $$
  so
  $$
    \Tra_{\Omega} [c_k] \otimes \zeta_m^{\otimes k} \equiv
    \Tra_{\Omega} [D_j^{(k)} y] \otimes (\zeta_m^j)^{\otimes k}
    \bmod P_{m,L/K}^{(k-1)} \otimes \mu_{p^m}^{\otimes k}.
  $$
  The first statement now follows by Theorem \ref{masseytrans}.

  Finally, we remark that
  replacing $a$ by $a^j$ but fixing $\sigma$ amounts to replacing
  $\zeta$ by $\zeta^j$, hence merely multiplies the Massey product
  by $j^k$, not affecting the group $P_{m,L/K}^{(k)}$.
\end{proof}

We now list the relationships between Massey products that result
when one of $m$, $n$, $K$, and $\Omega$ is allowed to vary.  We
leave the proofs to the reader.

\begin{lemma} \label{compat}
  Let
  $$
    Q_{m,\Omega/L/K}^{(k)} =
    H^2(G_{\Omega/K},\mu_{p^m}^{\otimes k+1})/P_{m,\Omega/L/K}^{(k-1)},
  $$
  and take $b \in U_{m,\Omega/L/K}^{(k)}$.
  \begin{enumerate}
    \item[a.]
      Fix $n' \le n$  with $m+r_k \le n'$, and let $L'$ be the unique
      subextension of $L/K$ of degree $p^{n'}$.  Then the quotient map
      $
        Q_{m,\Omega/L/K}^{(k)} \to Q_{m,\Omega/L'/K}^{(k)}
      $
      takes $(a,b)_{p^m,\Omega/L/K}^{(k)}$ to $(a,b)_{p^m,\Omega/L'/K}^{(k)}$.
    \item[b.]
      For $m' \le m$, the canonical map
      $
        Q_{m,\Omega/L/K}^{(k)} \to Q_{m',\Omega/L/K}^{(k)}
      $
      takes $(a,b)^{(k)}_{p^m,\Omega/K}$ to
      $(a,b)^{(k)}_{p^{m'},\Omega/K}$.
    \item[c.]
      Assume that $K'$ is an extension of $K$ contained in $\Omega$ with
      $L \cap K' = K$ and that $b'$ is an element of
      $U_{m,\Omega/LK'/K'}^{(k)}$ satisfying $N_{K'/K} b' = b$.  Then
      corestriction induces a map
      $
        Q_{m,\Omega/LK'/K'}^{(k)} \to Q_{m,n,\Omega/L/K}^{(k)}
      $
      taking $(a,b')^{(k)}_{p^m,\Omega/K'}$ to $(a,b)^{(k)}_{p^m,\Omega/K}$.
    \item[d.]
      If $\Omega'$ is a Galois extension of $K$ containing
      $\Omega$, then inflation induces a map
      $
        Q_{m,\Omega/L/K}^{(k)} \to Q_{m,n,\Omega'/L/K}^{(k)}
      $
      taking $(a,b)^{(k)}_{p^m,\Omega/K}$ to
      $(a,b)^{(k)}_{p^m,\Omega'/K}$.
  \end{enumerate}
\end{lemma}

\begin{remark}
  We have assumed that $[L:K] = p^n$.
  It would be nice to have the case of general
  degree written up for completeness.  We have not done this, as it
  would have complicated our discussion of Massey products, and it
  is not needed later in this paper.  Note that if
  $(a,b)^{(k)}_{p^m,\Omega/K}$ has a proper defining system, then
  one should have
  $$
    (a^i,b)^{(k)}_{p^m,\Omega/K} = i^k \cdot
    (a,b)^{(k)}_{p^m,\Omega/K} \pmod{ P_{m,L_i/K}^{(k-1)}},
  $$
  where $L_i = K(\alpha^i)$,
  for any $i \in \Z$ (not just prime to $p$).
\end{remark}

\section{Class groups} \label{classgroups}

We now take $K$ to be a number field, and we let $S$ be a set of
primes of $K$ including those above $p$ and all real places.
Recall that we assume $\mu_{p^n} \subset K^{\times}$.  We let
$G_{K,S}$ denote the Galois group of the maximal extension
$\Omega$ of $K$ unramified outside $S$, and we denote the Massey
products for $\Omega/K$ by $(a,b)_{p^m,K,S}^{(k)}$.

Let $\cl_{K,S}$ denote the $S$-class group of $K$, i.e., the
quotient of the class group of $K$ modulo the classes of all
finite primes in $S$.  Equivalently, it is the class group of the
ring of $S$-integers $\mc{O}_{K,S}$ of $K$, obtained from the integer ring of $K$ 
by inverting all finite primes in $S$. 
The cohomology group
$H^1(G_{K,S},\mathcal{O}_{\Omega,S}^{\times})$ is generated by the
classes of the cocycles $f_x$ for $x \in \Omega^{\times}$ with $x
\mathcal{O}_{\Omega,S} = \mathfrak{A} \mathcal{O}_{\Omega,S}$,
where $f_x(\tau) = x^{\tau-1}$ for $\tau \in G_{K,S}$ and
$\mathfrak{A}$ is a fractional ideal of $\OS{K}$. The map taking
the class of $f_x$ to the ideal class of $\mathfrak{A}$ provides a
canonical isomorphism
$H^1(G_{K,S},\mathcal{O}_{\Omega,S}^{\times}) \cong \cl_{K,S}$
(see, for example, \cite[Proposition 8.3.10]{nsw}). Via the
coboundary map
$$
  \partial \colon H^1(G_{K,S},\US{\Omega}) \to
  H^2(G_{K,S},\mu_{p^n}),
$$
we have a canonical injection
$$
  \Scl{K}/p^n\Scl{K} \to H^2(G_{K,S},\mu_{p^n}),
$$
and we denote the image of $\mf{A}$ in this group by $[\mf{A}]$.

\begin{lemma} \label{cotrans}
  Let $x \in \Omega^{\times}$ be such that $x \mathcal{O}_{\Omega,S}$ is the lift of a
  fractional ideal $\mathfrak{A}$ of $\mathcal{O}_{K,S}$.
  Then $\Tra_{\Omega} [x] = -[\mathfrak{A}].$
\end{lemma}

\begin{proof}
  Let $\xi^{p^n} = x$.
  Let $g$ denote a cochain with $p^n$th power equal to $f_x$, i.e., a cochain
  $g(\delta) = \xi^{\widetilde{\delta}-1}$, choosing a lift
  $\widetilde{\delta} \in G_K$ for each $\delta \in G_{K,S}$
  (with $\widetilde{1} = 1$).  By definition, we have that
  $$
     dg(\delta,\delta') =
     \frac{\xi^{\widetilde{\delta}-1}(\xi^{\widetilde{\delta'}-1})^{\widetilde{\delta}}}
     {\xi^{\widetilde{\delta\delta'}-1}}
     = \frac{\xi^{\widetilde{\delta}-1}\xi^{\widetilde{\delta}\widetilde{\delta'}-\widetilde{\delta}}}
     {\xi^{\widetilde{\delta\delta'}-1}}
     = \xi^{\widetilde{\delta}\widetilde{\delta'} - \widetilde{\delta\delta'}}.
  $$
  Define
  $$
    \theta(\delta,\delta') =
    \widetilde{\delta\delta'}^{-1}\widetilde{\delta}\widetilde{\delta'}.
  $$
  Letting $f_{\xi}$ denote the $p^n$th Kummer cocycle on $G_{\Omega}$ associated
  to $x$, we see that the cocycle defining $\mf{A}$ is given by
  $$
    dg(\delta,\delta') =
    (\xi^{\theta(\delta,\delta')-1})^{\widetilde{\delta\delta'}} =
    f_{\xi}(\theta(\delta,\delta')),
  $$
  since $\xi^{\theta(\delta,\delta')-1} \in \mu_{p^n}$.

  Now fix an extension $h \colon G_K \to \mu_{p^n}$ of $f_{\xi}$
  with $h(\widetilde{\delta}) = 1$ for all $\delta \in G_{K,S}$ and
  $$
    h(\widetilde{\delta}\tau) = h(\tau\widetilde{\delta}) = h(\tau)
  $$
  for $\tau \in G_K$.
  The coboundary $dh$ defines a cocycle on the quotient $G_{K,S}$,
  and the transgression $\Tra_{\Omega}[x]$ is given, since $h$ satisfies
  \eq{transprop}, by its class.  Since
  $$
    dh(\delta,\delta') = -h(\widetilde{\delta}\widetilde{\delta'})
    = -h(\widetilde{\delta\delta'}\theta(\delta,\delta'))
    = -h(\theta(\delta,\delta'))
    = -f_{\xi}(\theta(\delta,\delta')),
  $$
  we have our result.
\end{proof}

For a finite extension $F$ of $K$, we let $H_{F,S}$ denote the
$p$-Hilbert $S$-class field of $F$, which is the maximal
unramified abelian $p$-extension of $F$ in which all primes of $F$
above those in $S$ split completely.  Let $A_{F,S}$ denote the
$p$-part of $\cl_{F,S}$.  We have an identification $A_{F,S} \cong
\Gal(H_{F,S}/F)$ via class field theory.  We begin by recalling
the following basic lemma \cite[Lemma 6.1]{mcs}.

\begin{lemma}\label{genustheory}
  Let $L/K$ be a cyclic $p$-extension unramified outside $S$.
  The norm map $\Norm : A_{L,S} \to A_{K,S}$
  has image $\Gal(H_{K,S}/L \cap H_{K,S})$ and kernel
  equal to the subgroup of $A_{L,S}$ generated by $(\sigma-1)A_{L,S}$ and the
  intersection of $A_{L,S}$ with the subgroup of $\Gal(H_{L,S}/K)$ generated by
  all decomposition groups of primes in $S$.
\end{lemma}

Lemma \ref{genustheory} has the following easy corollary.

\begin{corollary} \label{genuscor}
  The map
  \begin{equation*}
    A_{L,S}/(\sigma-1)A_{L,S} \to A_{K,S},
  \end{equation*}
  induced by $\Norm$ is surjective if
  $L \cap H_{K,S} = K$ and injective if there is at most one prime
  in $S$ that does not split completely in $L/K$.
\end{corollary}

As in Section \ref{massey}, let $a \in K^{\times} \cap
\Omega^{\times p^n}$, $\alpha^{p^n} = a$, and $L = K(\alpha)$, and
assume $[L:K] = p^n$. Let $\sigma \in G_{L/K}$ be such that
$\sigma\alpha = \zeta\alpha$ for our fixed $p^n$th root of unity
$\zeta$.  For $i$ with $0 \le i \le n$, we let $L_i/K$ denote the
unique subextension of $L/K$ of degree $p^i$. Let $I_{F,S}$ denote
the $S$-ideal group of a finite extension $F$ of $K$.

\begin{theorem} \label{pairformula}
  Let $k$ and $m$ be positive integers with $k \le p^{n-m}(p-1)$.
  Let $b \in U_{m,L/K}^{(k)}$, and let $y \in L^{\times}$ with
  $N_{L/K} y = b$ and $D^{(k-1)}y \in \Omega^{\times p^m}$.
  Then we may write
  $$
    y\OS{L} \equiv \mf{Y}^{(\sigma-1)^k}\mf{B} \bmod p^mI_{L,S},
  $$
  with $\mf{Y} \in I_{L,S}$ and $\mf{B} \in I_{L_{n-m},S}$,
  and we have
  $$
    (a,b)^{(k)}_{p^m,K,S} \otimes \zeta_m^{\otimes -k} =
    -[N_{L/K} \mf{Y}] \pmod{P_{m,L/K}^{(k-1)}}
  $$
  unless $k = p^{n-m}(p-1)$, in which case
  $$
    (a,b)^{(k)}_{p^m,K,S} \otimes \zeta_m^{\otimes -k}
    = -[N_{L/K} \mf{Y}] - p^{m-1}[N_{L_{n-m}/K}\mf{B}]
    \pmod{P_{m,L/K}^{(k-1)}}.
  $$
\end{theorem}

\begin{proof}
  The first statement is a consequence of Lemma \ref{projform}, which
  can be seen as follows.  For each prime ideal $\mf{q}$ of
  $\mc{O}_{K,S}$, the $\Z[G_{L/K}]$-module summand $M_{\mf{q}}$
  of $I_{L,S}$ that is generated by the primes of
  $\mc{O}_{L,S}$ lying over $\mf{q}$ is cyclic, isomorphic
  to $\Z[H]$ for some quotient $H$ of $G_{L/K}$.
  Consider the projection $y_{\mf{q}}$ of $y\mc{O}_{L,S}$ to $M_{\mf{q}}$.
  Since $D^{(k-1)}y \in \Omega^{\times p^m}$, we have that
  $$
    D^{(k-1)}y_{\mf{q}} \equiv 0 \bmod p^mM_{\mf{q}}.
  $$
  We are now in the setting of Lemma \ref{projform}, with $R = \Z/p^m\Z$,
  $|H| = p^{n-s}$,
  and $G = G_{L/K}$.
  Identifying $R[H]$ with $M_{\mf{q}}/p^mM_{\mf{q}}$, we have that
  the image of $y_{\mf{q}}$ is contained in $\phi((I_G^k)^{\perp})^{\perp}$
  and therefore, we may write
  $$
    y_{\mf{q}} \equiv \mf{Y}_{\mf{q}}^{(\sigma-1)^k} \mf{B}_{\mf{q}}
    \mod p^mM_{\mf{q}}
  $$
  for some $\mf{Y}_{\mf{q}} \in M_{\mf{q}}$ and $\mf{B}_{\mf{q}} \in
  M_{\mf{q}}^{G_{L/L_{n-m}}}$.  Since $\mf{q}$ does not ramify in $L/K$,
  we have that
  $$
    M_{\mf{q}}^{G_{L/L_{n-m}}} = M_{\mf{q}} \cap I_{L_{n-m},S}.
  $$
  The claim then results from taking the product over all $\mf{q}$.

  As for the second statement, Lemma \ref{Dk} yields that
  $$
    D^{(k)}\mf{Y}^{(\sigma-1)^k} \equiv N_{L/K} \mf{Y} \bmod
    p^mI_{L,S}
  $$
  Furthermore, Lemma \ref{trivimage} implies that
  $D^{(k)}\mf{B} \in p^mI_{L,S}$ unless $k = p^{n-m}(p-1)$,
  in which case
  $$
    D^{(k)}\mf{B} \equiv p^{m-1}N_{L_{n-m}/K}\mf{B} \bmod
    p^mI_{L,S}.
  $$
  Recalling the first statement of the theorem
  and the fact that $\Omega$ contains the $p$-Hilbert class field
  of $L$ (and without loss of generality, replacing $b$ by a prime-to-$p$
  power of itself), we conclude that
  there exists $w \in \Omega^{\times}$ with
  \begin{equation} \label{ideal}
    (w^{p^m}c_k) \mathcal{O}_{\Omega,S} =
    \begin{cases}
      N_{L/K} \mf{Y} & \mr{if\ } k < p^{n-m}(p-1) \\
      N_{L/K} \mf{Y} \cdot N_{L_{n-m}/K} \mf{B}^{p^{m-1}} & \mr{if\ } k
      = p^{n-m}(p-1),
    \end{cases}
  \end{equation}
  for $c_k = D^{(k)}y$ as before.
  By Theorem \ref{masseytrans}, we need only
  make the following computation in cohomology with
  $\mu_{p^m}$-coefficients:
  \begin{equation} \label{eqcocyc}
    \Tra_{\Omega} [c_k] = \Tra_{\Omega} [w^{p^m} c_k]
    =
    \begin{cases}
      -[N_{L/K} \mf{Y}] & \mr{if\ } k < p^{n-m}(p-1) \\
      -[N_{L/K} \mf{Y}] - p^{m-1}[N_{L_{n-m}/K} \mf{B}] & \mr{if\ } k = p^{n-m}(p-1),
    \end{cases}
  \end{equation}
  where the last step follows from \eq{ideal} and
  Lemma \ref{cotrans}.
\end{proof}

\begin{remark}
  The case $k = 1$ and $m = n$ of the above result is Theorem 2.4 of
  \cite{mcs} on cup products, giving an entirely different proof of
  that result which, when boiled down to what is needed for that
  specific case, seems simpler.  (Actually, we have restricted to
  the case that the degree is $p^n$ here, but for cup products, the
  general case follows almost immediately.)
\end{remark}

For $m \le n$ and $k \ge 1$, let
\begin{equation} \label{Jk}
  J_{m,L/K}^{(k)} = \langle \mf{a} \in A_{L,S} \mid
  \mf{a}^{(\sigma-1)^k} \in \phi(A_{L_{n-m},S})
  + p^m A_{L,S} \rangle,
\end{equation}
where $\phi \colon A_{L_{n-m},S} \to A_{L,S}$ is the natural map.

\begin{proposition} \label{isogroups}
  Let $k \ge 1$ and $m \ge 1$ be such that $k \le p^{n-m}(p-1)-1$.
  The map $N_{L/K}: A_{L,S} \to A_{K,S}$ induces a
  surjective map
  $$
    J_{m,L/K}^{(k)}/(J_{m,L/K}^{(k)} \cap (p^m,\sigma-1)A_{L,S})
    \to P_{m,L/K}^{(k)}
  $$
  which is an isomorphism if there is a prime $v$ of $K$ in $S$
  that does not split at all in $L$ and $v$ is the only prime in $S$
  which does not split completely in $L$.
\end{proposition}

\begin{proof}
  Given $\mf{a} \in I_{L,S}$ with ideal class in $J_{m,L/K}^{(k)}$, we
  can find $b \in K^{\times}$ associated with $\mf{a}$ as follows.
  Apply $(\sigma-1)^k$ to $\mf{a}$, multiply it by an element of
  $I_{L_{n-m},S} \cdot I_{L,S}^{p^m}$
  to make it principal,
  and apply $N_{L/K}$ to a generator of the resulting principal ideal
  to obtain $b$.
  One sees that $b \in U_{m,L/K}^{(k)}$, and then
  Theorem~\ref{pairformula} implies that $[N_{L/K}\mf{a}] \in P_{m,L/K}^{(k)}$.
  Hence, $N_{L/K} J_{m,L/K}^{(k)}$ modulo $p^m$ is contained in
  $P_{m,L/K}^{(k)}$.

  Conversely, take an element $t \in P_{m,L/K}^{(k)}$ associated to a proper
  defining system for some $b \in U_{m,L/K}^{(k)}$ and $y$ as in
  Theorem \ref{pairformula}.  Then the first statement of Theorem
  \ref{pairformula} supplies an ideal $\mathfrak{a}$ with class in
  $J_{m,L/K}^{(k)}$ such that $t = [\Norm \mathfrak{a}]$, noting the equality of
  cohomology classes demonstrated in \eq{eqcocyc}.  Therefore, we have
  surjectivity.

  Finally, if $v$ is as stated, then the map in
  Corollary~\ref{genuscor} is an isomorphism, so the kernel of
  the map $A_{L,S} \to A_{K,S}/p^m$ induced by norm is given by
  $(p^m,\sigma-1)A_{L,S}$.  Thus we have injectivity.
\end{proof}

\section{Iwasawa theory} \label{infiniteext}

Now let $F$ be a number field, and let $K = F(\mu_p)$ if $p$ is odd
and $K = F(\mu_4)$ if $p = 2$.  
In this section, we pass to the cyclotomic $\zp$-extension $K_{\infty} =
K(\mu_{p^{\infty}})$ of $K$, and we consider a $\zp$-extension
$L_{\infty}$ of $K_{\infty}$ unramified outside a finite set of
primes of $K_{\infty}$.  We will consider Massey products in the
inverse limits of class groups up the cyclotomic tower.

Fix any set of primes $S$ of $K$ containing those
above $p$ and those which ramify in $L_{\infty}/K$.
Let $\zeta = (\zeta_n)$ be a generator of the Tate module, and let
$K_n = K(\mu_{p^n})$ for $n \ge 1$. There exists a sequence of
elements $a = (a_n)_{n \ge N}$ for $N$ sufficiently large and a
nondecreasing sequence $(l_n)_{n \ge N}$ with infinite limit and
$l_n \le n$ for all $n$, such that $a_n \in K_n^{\times} \cap
\Omega^{\times p^{l_n}}$, $a_n \notin K_n^{\times p}$, with
$$
  a_{n+1}a_n^{-1} \in K_{n+1}^{\times p^{l_n}}
$$
and such that, setting $\alpha_n^{p^{l_n}} = a_n$ and $L_n =
K_n(\alpha_n)$, the field $L_{\infty}$ is the union of the $L_n$.

\begin{remark}
    Any prime $w$ in $K_{\infty}$ lying over a prime
    other than $p$ which is unramified in a pro-$p$ extension $F/E$ with
    $E$ containing $K_{\infty}$ must split completely in $F/E$, since
    the completion of $K_{\infty}$ at $w$ contains the unramified
    $\zp$-extension of the completion $K_w$.
\end{remark}

Given this remark, we let $\mc{A}_{K}$ denote the Galois group of the maximal
pro-$p$ unramified abelian extension of $K_{\infty}$ in which all primes
above $p$ split completely, and let $\mc{A}_{L}$ be similarly
defined for $L_{\infty}$.
Class field theory provides canonical
isomorphisms
$$
  \mc{A}_{K} \cong \lim_{\leftarrow} A_{K_n,S} \cong
  \lim_{\leftarrow} A_{K_n,S}/p^{m_n},
$$
with inverse limits taken with respect to norm maps, for any
nondecreasing sequence $m_n$ with infinite limit.  Similarly,
$\mc{A}_{L}$ may be viewed as an inverse limit of the
$A_{L_n,S}$ (modulo $p^{m_n}$) via norm maps, even though the
extensions may not be Galois.

Let $G = G_{L_{\infty}/K_{\infty}}$, and let $I_G$ denote the
augmentation ideal of $\Z_p[[G]]$.   
Let $\tilde{\Lambda} = \zp[[G_{K_{\infty}/F}]]$.  We remark that
restriction induces a homomorphism
$$
  \mc{A}_{L}/I_G\mc{A}_{L} \to \mc{A}_{K}
$$
of $\Z_p$-modules which is a homomorphism of $\tilde{\Lambda}$-modules if
$L_{\infty}/K$ is Galois.

Let $m_n$ be any (fixed) nondecreasing sequence with infinite
limit and satisfying $k < p^{l_n-m_n}(p-1)$ for all (sufficiently
large) $n$.  We define
$$
  \mc{P}_{L/K}^{(k)} = \lim_{\leftarrow} P_{m_n,L_n/K_n}^{(k)}.
$$
That $\mc{P}_{L/K}^{(k)}$ is well-defined and its definition is
independent of the choices of $l_n$ and $m_n$ follows from the
first three parts of Lemma \ref{compat}.  It is also independent
of the choice of $S$, as the following lemma will make clear.
Define $\mc{J}_{L/K}^{(k)}$ to be the kernel of $(\sigma-1)^k$ on
$\mc{A}_{L}$, where $\sigma$ is a topological generator of $G$.

\begin{lemma} \label{infinitecase}
  The restriction map $\mc{A}_{L} \to \mc{A}_{K}$
  induces a surjective homomorphism of $\Z_p$-modules
  $($of $\tilde{\Lambda}$-modules if $L_{\infty}/K$ is Galois$)$:
  \begin{equation} \label{infJk}
    \mc{J}_{L/K}^{(k)}/(\mc{J}_{L/K}^{(k)} \cap
    I_G\mc{A}_{L}) \to \mc{P}_{L/K}^{(k)}.
  \end{equation}
\end{lemma}

\begin{proof}
  We begin by verifying that
  \begin{equation} \label{inverselim}
    \mc{J}_{L/K}^{(k)}/(\mc{J}_{L/K}^{(k)} \cap I_G\mc{A}_{L}) \cong
    \lim_{\leftarrow} J_n/(J_n \cap (p^{m_n} + I_G) A_{L_n,S}),
  \end{equation}
  where we define $J_n$ to be $J_{m_n,L_n/K_n}^{(k)}$ as defined in \eq{Jk}.
  The inverse limit in \eq{inverselim}
  must be interpreted appropriately as sequences
  of norm compatible elements of the subgroups $J_n$ of $A_{L_n,S}$, since
  $N_{L_t/L_n} J_t$ need not be contained in $J_n$ for $t > n$.

  By definition, the image of $\mc{J}_{L/K}^{(k)}$ under the natural map to $A_{L_n,S}$
  is contained in $J_n$.  Conversely,
  take a sequence $\mf{a} = (\mf{a}_n)$ of norm compatible elements
  $$
    \mf{a}_n \in J_n + (p^{m_n}+I_G) A_{L_n,S}.
  $$
  Let $t > n$, and let
  $\phi_t$ denote the natural map from the $S$-class group of an extension
  of $K$ in $L_t$ to that of $L_t$.
  We have by definition of $J_t$ that
  $$
    \mf{a}_n^{(\sigma-1)^k} = N_{L_t/L_n} (\mf{a}_t^{(\sigma-1)^k})
    \in N_{L_t/L_n} \phi_t(A_{L_{t-m_t},S})
    + p^{m_t} A_{L_n,S} + I_G^{k+1} A_{L_n,S}.
  $$
  Note that if $t$ is such that $t \ge m_t+n$, then we have
  $$
    N_{L_t/L_n} \phi_t(A_{L_{t-m_t},S}) = p^{m_t} N_{L_{t-m_t}/L_n} A_{L_{t-m_t},S}
    \subseteq p^{m_t} A_{L_n,S}.
  $$
  Therefore, noting that $m_t$ has infinite limit, we have
  $$
    \mf{a}_n^{(\sigma-1)^k} \in \bigcap_t \left( p^{m_t} A_{L_n,S} + I_G^{k+1}
    A_{L_n,S} \right) = I_G^{k+1} A_{L_n,S}.
  $$
  This yields that
  $\mf{a}^{(\sigma-1)^k} \in I_G^{k+1} \mc{A}_{L}$,
  and therefore, that
  $\mf{a} \in \mc{J}_{L/K}^{(k)} + I_G \mc{A}_{L},$
  finishing the verification of \eq{inverselim}.

  Now consider the inverse limit of the surjective maps
  \begin{equation*} \label{surjmaps}
    \psi_n \colon J_n/(J_n \cap (p^{m_n}+I_G)A_{L_n,S}) \to
    P_{m_n,L_n/K_n}^{(k)}
  \end{equation*}
  of Proposition \ref{isogroups} (for $n$ sufficiently large).
  The inverse limit of the $\psi_n$ is also surjective, as
  the kernels of the $\psi_n$ are finite abelian groups.
  The map in \eq{infJk} is a surjection as
  the composition of the isomorphism in \eq{inverselim} with
  $\displaystyle \lim_{\leftarrow} \psi_n$.
\end{proof}

Let $\mc{H}_{K}$ and $\mc{H}_{L}$ denote the fields such that
\begin{eqnarray*}
\mc{A}_{K} = G_{\mc{H}_{K}/K_{\infty}} & \mr{and} &
\mc{A}_{L} = G_{\mc{H}_{L}/L_{\infty}}.
\end{eqnarray*}
Let $\mc{D}_{L/K}$ be the intersection of $\mc{A}_{L}$ with the
subgroup of $G_{\mc{H}_{L}/K_{\infty}}$ generated by all decomposition
groups of primes of $K_{\infty}$.
We also define
$$
  \mc{A}_{L/K} = \Gal(\mc{H}_{K}/\mc{H}_{K} \cap L_{\infty}) \le
  \mc{A}_{K}.
$$

\begin{lemma} \label{genuszp}
  The restriction map
  \begin{equation*} \label{normclass}
    \mc{A}_{L}/I_G\mc{A}_{L} \to \mc{A}_{K}
  \end{equation*}
  has kernel $\mc{D}_{L/K}/(\mc{D}_{L/K} \cap I_G\mc{A}_{L})$
  and image $\mc{A}_{L/K}$.
\end{lemma}

\begin{proof}
  The result follows immediately from Lemma \ref{genustheory} by
  taking inverse limits.
\end{proof}

Note that Lemma \ref{infinitecase} implies that $\mc{A}_{L/K}$
contains $\mc{P}_{L/K}^{(k)}$ for all $k \ge 0$. We define the
{\em decomposition-free quotient} of the $k$th graded piece 
$$
	\mr{gr}_G^k \mc{A}_L = I_G^k \mc{A}_L/I_G^{k+1} \mc{A}_L 
$$
in the augmentation filtration of $\mc{A}_L$ by
$$
   \mc{Q}_{L/K}^{(k)} =
   I_G^k \mc{A}_{L}/(I_G^k \mc{D}_{L/K}+I_G^{k+1} \mc{A}_{L})
$$
for each $k \ge 0$.  Regarding these quotients, we have the following
obvious remark (which follows from Corollary \ref{genuscor}).

\begin{remark}
  If at most one prime above $p$ in $K_{\infty}$ does not split completely
  in $L_{\infty}$, then in fact
  $\mc{Q}_{L/K}^{(k)} = \mr{gr}_G^k \mc{A}_{L}$ for all $k$.
\end{remark}

Our main result now describes the decomposition-free quotients in
terms of Massey products.

\begin{theorem} \label{mainthm}
  For any $k \ge 0$, we have a canonical isomorphism of $\Z_p$-modules
  $$
    \mc{Q}_{L/K}^{(k)}
    \xrightarrow{\sim} \mc{A}_{L/K}/\mc{P}_{L/K}^{(k)} \otimes_{\Z_p}
    G^{\otimes k}
  $$
  which is an isomorphism of $\tilde{\Lambda}$-modules if $L_{\infty}/K$ is Galois.
\end{theorem}

\begin{proof}
  Taking the quotient of the surjection
  $$
    \mc{A}_{L}/I_G\mc{A}_{L} \to \mc{A}_{L/K}
  $$
  by the surjection in Lemma \ref{infinitecase}
  yields a surjection
  $$
    \mc{A}_{L}/(\mc{J}_{L/K}^{(k)} + I_G\mc{A}_{L}) \to
    \mc{A}_{L/K}/\mc{P}_{L/K}^{(k)}
  $$
  with kernel given by the image of $\mc{D}_{L/K}$.
  By definition of $\mc{Q}_{L/K}^{(k)}$, we need only observe
  that there is a canonical isomorphism
  $$
    \mc{A}_{L}/(\mc{J}_{L/K}^{(k)} + I_G\mc{A}_{L})
    \otimes_{\zp} G^{\otimes k}
    \xrightarrow{\sim} \mr{gr}_G^k\mc{A}_{L}  
  $$
  given by
  $
    \mf{a} \otimes \sigma^{\otimes k} \mapsto
    \mf{a}^{(\sigma-1)^k}.
  $
\end{proof}

This has the following immediate corollary, describing a setting in which we
have an exact comparison between the augmentation filtration of $\mc{A}_{L}$
and the group $\mc{A}_K$ modulo the image of Massey products.

\begin{corollary}
  Assume that there exists a unique prime $v$ of $K_{\infty}$ that does not
  split completely in $L_{\infty}$ and that $v$ does not split
  at all in $L_{\infty}$.  Then there are canonical isomorphisms
  $$
    \mr{gr}_G^k \mc{A}_L \xrightarrow{\sim}
    \mc{A}_K/\mc{P}_{L/K}^{(k)} \otimes_{\zp} G^{\otimes k}
  $$
  for all $k \ge 0$.
\end{corollary}

Let
$$
  \mc{U}_K = \lim_{\leftarrow} \, (\US{K_n} \otimes \zp),
$$
the inverse limit being taken with respect to norm maps.  This
group is easily seen to satisfy
$$
  \mc{U}_K \cong \lim_{\leftarrow} \, \US{K_n}/\OS{K_n}^{\times p^n} \cong
  \lim_{\leftarrow} \, K_n^{\times}/K_n^{\times p^n},
$$
the latter isomorphism following from the fact that no prime
splits completely in the cyclotomic $\zp$-extension a number field.  We
write an element of $\mc{U}_K$ as a sequence $b = (b_n)$ of elements
$b_n \in \US{K_n}$ with reduction modulo $p^n$th powers forming a
norm compatible sequence.

For $k \ge 1$, let
$$
  \mc{U}_{L/K}^{(k)} = \langle b \in \mc{U}_K \mid b_n \in
  U_{m_n,L_n/K_n}^{(k)} \mr{\ for\ all\ } n \gg 0
  \rangle.
$$
For $b \in \mc{U}_{L/K}^{(k)}$, we define Massey products
$$\label{infinitemassey}
  (a,b)^{(k)}_{p^{\infty},K_{\infty}} \otimes \zeta^{\otimes -k} =
  \lim_{\leftarrow}\, (a_n,b_n)^{(k)}_{p^{m_n},K_n,S} \otimes
  \zeta_{m_n}^{\otimes -k} \in \mc{A}_{K}/\mc{P}_{L/K}^{(k-1)}.
$$
One would like to say that $\mc{P}_{L/K}^{(k)}$ is generated (modulo
$\mc{P}_{L/K}^{(k-1)}$) by such Massey products.  This is indeed the
case.

\begin{proposition} \label{pairgroup}
  For $k \ge 1$, we have an equality of groups
  $$
    \mc{P}_{L/K}^{(k)}/\mc{P}_{L/K}^{(k-1)} =
    \langle (a,b)^{(k)}_{p^{\infty},K_{\infty}}
    \otimes \zeta^{\otimes -k} \mid b \in \mc{U}_{L/K}^{(k)}
    \rangle.
  $$
\end{proposition}

\begin{proof}
  Let $\mf{A} \in \mc{P}_{L/K}^{(k)}$.
  Then $\mf{A}$ is the image of some
  $\mf{a} = (\mf{a}_n) \in \mc{J}_{L/K}^{(k)}$ by Lemma \ref{infinitecase}.
  Pick a representative of $\mf{a}_n$, apply $(\sigma-1)^k$ to it, and let $y_n$ be
  a generator.  The element $b_n = N_{L_n/K_n} y_n$ is determined by
  $\mf{a}$ up to the norm of an element of $\US{L_n}$.  The sequence
  $(b_n)$ therefore yields a norm compatible sequence of elements in
  $$
    \lim_{\leftarrow}\, U_{m_n,L_n/K_n}^{(k)}/N_{L_n/K_n}\US{L_n}.
  $$
  By Theorem \ref{pairformula}, we have that $\mf{A}
  \pmod{\mc{P}_{L/K}^{(k-1)}}$
  is the inverse limit of the Massey products
  $$
    (a_n,b_n)^{(k)}_{p^{m_n},K_n,S} \otimes \zeta_{m_n}^{\otimes -k}.
  $$
  We have only to check that the map
  $$
    \mc{U}_{L/K}^{(k)} \cong \lim_{\leftarrow}  \,
    U_{m_n,L_n/K_n}^{(k)}/\OS{K_n}^{\times p^{m_n}}
    \to \lim_{\leftarrow} \,
    U_{m_n,L_n/K_n}^{(k)}/(\OS{K_n}^{\times p^{m_n}} \cdot
    N_{L_n/K_n}\US{L_n})
  $$
  is surjective, but this follows as $U_{m_n,L_n/K_n}^{(k)}$ may be
  replaced by $U_{m_n,L_n/K_n}^{(k)} \cap \US{K_n}$ in the inverse limits 
  and $\US{K_n}/\OS{K_n}^{\times p^{m_n}}$ is finite.
\end{proof}

We conclude this section with the following statement on which
elements of $\mc{U}$ are trivial under all Massey products with
$a$.

\begin{proposition} \label{trivmassey}
  Let $\mc{U}_{L/K} = \cap_{k \ge 1} \mc{U}_{L/K}^{(k)} \subseteq \mc{U}_K$.  Then
  \begin{equation} \label{sameunitgp}
    \mc{U}_{L/K} =
    \langle b \in \mc{U} \mid b_n \in N_{L_n/K_n} (\US{L_n}
    \otimes \zp)
    \mr{\ for\ all\ } n \ge N \rangle.
  \end{equation}
\end{proposition}

\begin{proof}
  The group on the right-hand side of
  \eq{sameunitgp} is contained in
  $\mc{U}_{L/K}$ by Theorem \ref{pairformula}.
  Assume,
  on the other hand, that $b \in \mc{U}_{L/K}$.
  Fix $n \ge N$ and choose $k$ large enough so that $I_G^{k-1}A_{L_n,S}
  = 0$.  For any $t \ge n$ satisfying $t-m_t \ge r_k$, we have that
  $b_t \in U_{m_t,L_t/K_t}^{(k)}$.  Hence, Theorem \ref{pairformula}
  yields that $b_t = N_{L/K} y_t$
  for some $y_t \in L_t^{\times}$ with $y_t\OS{L_t} =
  \mf{Y}_t^{(\sigma-1)^k}\mf{B}_t$, where $\mf{Y}_t \in I_{L_t,S}$ and
  $\mf{B}_t \in I_{L_{t-m_t},S}$.
  Let $y_n = N_{L_t/L_n} y_t$.  Then we have
  $$
    y_n \mc{O}_{L_n,S} = (N_{L_t/L_n} \mf{Y}_t)^{(\sigma-1)^k}
    N_{L_t/L_n} \mf{B}_t
  $$
  and $N_{L_n/K_n} y_n = b_n v^{p^n}$ for some $v \in
  \mc{O}_{K_n,S}^{\times}$.
  If $t-m_t \ge n$ as well, then
  $N_{L_t/L_n} \mf{B}_t \in p^{m_t}I_{L_n,S}.$
  As in the proof of Lemma \ref{infinitecase}, since $m_t$ can be made arbitrarily
  large and $n$ was arbitrary, we see that
  $y_n\mc{O}_{L_n,S} = \mf{C}^{(\sigma-1)^k}$
  for some $\mf{C} \in I_{L_n,S}$.
  Since $I_G^{k-1} A_{L_n,S}$ is trivial by
  assumption, some prime to $p$-power of $y_n$ (which we assume without
  loss of generality to be $1$) has the form
  $u \cdot x^{\sigma-1}$ for some $u \in \US{L_n}$
  and $x \in L_n^{\times}$, which tells us that $b_n =
  N_{L_n/K_n} (uv^{-1})$.
\end{proof}

\begin{remark}
  Let $\mc{U}_L$ be the group
  $\displaystyle \lim_{\leftarrow}\, (\US{L_n} \otimes \zp)$
  of universal norm sequences from $L_{\infty}$.
  There is an obvious norm map
  $$
    N_{L_{\infty}/K_{\infty}} \colon \mc{U}_L \to \mc{U}_K.
  $$
  As a consequence of Proposition \ref{trivmassey}, we have
  $$
    \mc{U}_{L/K} = N_{L_{\infty}/K_{\infty}} \mc{U}_L,
  $$
  the latter group being easily seen to be equal to the group
  on the right-hand side of \eq{sameunitgp}.
\end{remark}

\section{Examples} \label{examples}

In this section, we apply the results of Section \ref{infiniteext}
to the case in which $F = \Q$ and $K = \Q(\mu_p)$ for an odd prime $p$. 
Let $S$ be the set consisting of the unique prime
$(1-\zeta_1)$ above $p$ in $K$.  
We take $L_{\infty}$ to be a $\zp$-extension of $K_{\infty} =
\Q(\mu_{p^{\infty}})$ which is Galois over $\Q$ and unramified
outside $p$.  Let $G =
G_{L_{\infty}/K_{\infty}}$ and $\mc{G} = G_{L_{\infty}/\Q}$. We
remark that $\mc{A}_{K}$ is equal to the Galois group of the
maximal pro-$p$ unramified abelian extension of $K_{\infty}$.

Let $T = \zp \times \Z/(p-1)\Z$.  We view $\Z$ as a subgroup of
$T$ via the diagonal map.  For $n \ge 1$, we have a natural
surjection $T \to \Z/p^{n-1}(p-1)\Z$ which allows us to speak of
congruences modulo $p^{n-1}(p-1)$ between elements of $T$, and we
let $t_n$ denote a lift to $\Z$ of the 
image of $t \in T$ under this map.  For $t
\in T$, define
$$
  \zp(t) = \lim_{\leftarrow}\, (\Z/p^n\Z)(t_n).
$$

Let $\Delta = \Gal(K/\Q)$, let $\omega$ denote the Teichm\"uller
character, and set
$$
  \epsilon_i = \sum_{\delta \in \Delta} \omega(\delta)^{-i}\delta
  \in \zp[\Delta]
$$
for $i \in T$.  We
have the following corollary of Theorem \ref{mainthm}.

\begin{lemma} \label{consequence}
  For $i \in \Z$ and $k \ge 0$, there is a canonical isomorphism of
  $\tilde{\Lambda}$-modules:
  $$
    \epsilon_{i+kt}(\mr{gr}_G^k \mc{A}_{L})
    \otimes_{\Z_p} G^{\otimes -k} \cong
    \epsilon_i(\mc{A}_{L/K}/\mc{P}_{L/K}^{(k)})
  $$
\end{lemma}
 
For $t \in T$ with $t$ odd (modulo
$p-1$), define
$$
  \lambda_{t,n} = \prod_{\substack{i = 1\\(i,p)=1}}^{p^n-1}
  (1-\zeta_n^i)^{i^{t_n-1}}.
$$
Then
$$
  \lambda_{t,n+1}\lambda_{t,n}^{-1} \in K_{n+1}^{\times p^n}.
$$
We fix an odd $t \in T$ and assume for now that $L_{\infty}/K_{\infty}$ is
defined by the sequence $\lambda_t = (\lambda_{t,n})_{n \ge 1}$.
The field $L_{\infty}$ is then a $\zp$-extension of $K_{\infty}$
unramified outside $p$ which satisfies $G \cong \zp(t)$ as a
$\tilde{\Lambda}$-module.  We remark that
$L_{\infty}/K_{\infty}$ is totally ramified if $\epsilon_t A_K =
0$.

Let $A_K$ be the $p$-part of the class group of $K$, so $A_K = A_{K,S}$.
Recall that Vandiver's conjecture states that the real part
$A_K^+$ of $A_K$ is zero and implies the weaker statement that
$\epsilon_s A_K$ is cyclic for all odd $s$.  Vandiver's conjecture
is known to hold for all $p <$ 12,000,000 \cite{bcems}.  We assume the
weaker condition of
cyclicity of an odd eigenspace of $A_K$ in the following
proposition (see \cite{kurihara} for a discussion of this
condition).  For $r \in \Z$, $n \ge 1$ (or $\infty$), and $a, b
\in \mc{O}_{K_n,S}^{\times}$, we let
$$
  \langle a,b \rangle_{p^n,r} = \epsilon_{2-r}(a,b)^{(1)}_{p^n,K_n,S}
$$
(see \cite{mcs} for a study of this cup product pairing in the
case $n = 1$).

\begin{proposition} \label{smallclassgp}
  Let $r$ be even such
  that $\epsilon_{1-r}A_K$ has $p$-rank $1$.  Let $M$ be a
  $\zp[[\mc{G}]]$-submodule of
  $\mc{A}_{L}$ generated by a choice of element with restriction
  generating $\epsilon_{1-r}\mc{A}_{L/K}$.  If
  $\langle \lambda_{t,1}, \lambda_{r-t,1} \rangle_{p,r}$ is nonzero,
  then $M \cong \epsilon_{1-r}\mc{A}_K$.
\end{proposition}

\begin{proof}
  The condition of $\epsilon_{1-r}A_K$ having
  $p$-rank $1$ is equivalent to that of
  $\epsilon_{1-r}\mc{A}_K$ being principally generated as a
  $\Lambda = \zp[[G_{K_{\infty}/K}]]$-module \cite[Remark 1.5]{kurihara}.  Set
  $$
    X = \epsilon_{1-r}(\mc{A}_K/\mc{P}_{L/K}^{(1)}).
  $$
  Let $\mf{m}$ denote the maximal ideal of $\Lambda$.
  The assumption on the pairing implies that
  $$
    \epsilon_{2-r}(A_K \otimes \mu_p)/\langle\langle \lambda_{t,1},\lambda_{r-t,1}
    \rangle_{p,r}\rangle
  $$
  is trivial.  Now,
  $X/\mf{m}X$ is isomorphic to a quotient of this, hence is
  trivial, which in turn implies $X = 0$.
  On the other hand, Lemma \ref{consequence} yields
  $$
    \mr{gr}^1_G M \cong \epsilon_{1-r+t}(\mr{gr}^1_G \mc{A}_{L})
    \cong \epsilon_{1-r}(\mc{A}_{L/K}/\mc{P}_{L/K}^{(1)})
    \hookrightarrow X.
  $$
  Hence, we have $I_G M = 0$ and $\epsilon_{1-r}\mc{A}_{L/K} =
  \epsilon_{1-r}\mc{A}_K$,
  and, therefore, we have
  $$
    M = M/I_G M \cong \epsilon_{1-r}\mc{A}_K.
  $$
\end{proof}

We now apply Proposition \ref{smallclassgp} to one case in which
the pairing $\langle \, \cdot \, , \, \cdot \, \rangle_{p,r}$ has
been computed (see \cite{mcs}).  Let $\mc{Z}_L$ denote the Galois
group of the maximal unramified 
abelian pro-$p$ extension of $L_{\infty}$.

\begin{proposition} \label{p37}
  For $p = 37$ and $t \in T$ with $t \not\equiv 5,27 \bmod 36$,
  we have isomorphisms of $\zp[[\mc{G}]]$-modules,
  $$
    \mc{A}_{L} \cong \mc{A}_K \cong \zp(s),
  $$
  for some $s \in T$ with
  $s \equiv 5 \bmod 36$.  If, furthermore,
  $t \not\equiv 31 \bmod 36$, then $\mc{Z}_L \cong \zp(s)$
  as well.
\end{proposition}

\begin{proof}
  Note that $p \mid B_{32}$.  By \cite[Theorem 7.5]{mcs}, the
  pairing value
  $\langle \lambda_{t,1}, \lambda_{32-t,1} \rangle_{37,32}$ is
  nonzero unless $t \not\equiv 5,27 \bmod 36$.  (Note that if
  $t \equiv 5 \bmod 36$, then $L_{\infty}/K_{\infty}$ is not
  totally ramified at $37$.) It follows from
  Proposition \ref{smallclassgp} that $\mc{A}_{L} \cong
  \mc{A}_K$.  Furthermore, the $\lambda$-invariant of
  $\mc{A}_K = \epsilon_5 \mc{A}_K$ is $1$, so $\mc{A}_K \cong
  \zp(s)$ for some $s \equiv 5 \bmod 36$.

  Note that
  $$
    \mc{A}_{L} \cong \mc{Z}_L/\langle \mf{p} \rangle,
  $$
  where $\mf{p}$ denotes the class of the unique
  prime of $L_{\infty}$ above $p$.  Since $L_{\infty}/K_{\infty}$ is
  totally ramified at $p$, we have
  $$
    \mc{Z}_L/I_G\mc{Z}_L \cong \mc{A}_{L}/I_G\mc{A}_{L} \cong \mc{A}_K.
  $$
  In particular, we have $\mf{p} \in I_G\mc{Z}_L$.  The fact that
  $I_G\mc{A}_{L} = 0$ implies that $I_G\mc{Z}_L$
  is generated by $\mf{p}$ and that
  $$
    I_G\mc{Z}_L \cong
    \epsilon_{5+t}(I_G\mc{Z}_L/I_G^2\mc{Z}_L).
  $$
  Since $\mf{p} \in \mc{Z}_L^{\mc{G}}$, we conclude that $I_G \mc{Z}_L = 0$
  if $t \not\equiv 31 \bmod 36$.
\end{proof}

We now illustrate that $\mc{A}_{L}$ can have larger $p$-rank
than $\mc{A}_K$ as well.

\begin{proposition} \label{larger}
  Let $r \in T$, and assume that
  $\epsilon_{1-r}\mc{A}_K \cong \zp(1-r)$ and
  $\epsilon_r\mc{A}_K = \epsilon_{2r-1}\mc{A}_K = 0$.
  Let $L_{\infty}$ be the $\zp$-extension of $K_{\infty}$
  defined by $\lambda_{2r-1}$.
  Then we have
  $$
    \epsilon_{1-r}(\mr{gr}_G^k \mc{A}_L \otimes_{\zp}
    G^{\otimes -k}) \cong \zp(1-r).
  $$
  if $k = 0$ or $1$, and for $k = 2$ if $3r \neq 2$.
\end{proposition}

\begin{proof}
  For $k = 0$, the result is Lemma \ref{genuszp}.

  For $k = 1$, let
  $1-\zeta$ denote the universal norm for $K_{\infty}/K$
  given by the sequence $(1-\zeta_n)$.
  We remark that $\epsilon_{1-r}\mc{P}_{L/K}^{(1)}$
  is generated as a $\Lambda$-module by
  $\langle \lambda_{2r-1},1-\zeta \rangle_{p^{\infty},r}$, since
  $\epsilon_r \mc{A}_K = 0$ implies that $\epsilon_r\mc{U}_K$ is
  generated by universal norms arising from cyclotomic units.
  At the level of $K_n$, we have
  \begin{equation} \label{samepair}
    \langle \lambda_{2r-1,n},1-\zeta_n \rangle_{p^n,r} =
    \langle 1-\zeta_n,\lambda_{1-r,n} \rangle_{p^n,r},
  \end{equation}
  from which we conclude that
  $$
    \epsilon_{1-r}\mc{P}_{L/K}^{(1)} =
    \epsilon_{1-r}\mc{P}_{H/K}^{(1)},
  $$
  where $H_{\infty}/K_{\infty}$ is the $\zp$-extension defined by
  $\lambda_{1-r}$.
  The extension $H_{\infty}/K_{\infty}$ is unramified.  Using Lemma
  \ref{genustheory},
  and
  defining $\mc{A}_{H}$ to be the Galois group of the maximal pro-$p$ abelian
  unramified extension of $H_{\infty}$ in which all primes above $p$ split
  completely,
  we see by Lemma \ref{genuszp} that
  $$
    \epsilon_{1-r}(\mc{A}_{H}/I_G \mc{A}_{H}) = 0.
  $$
  Since the right hand term of
  \eq{samepair} generates the $\epsilon_{1-r}$-eigenspace of
  $P_{n,H_n/K_n}^{(1)}$, where $H_n/K_n$ is the extension given
  by a $p^n$th root of $\lambda_{1-r,n}$.  Since
  $$
    \mc{P}_{H/K}^{(1)} = \lim_{\leftarrow}
    P_{n,H_n/K_n}^{(1)},
  $$
  Theorem \ref{pairformula} implies that
  $\epsilon_{1-r} \mc{P}_{H/K}^{(1)}$ is trivial.

  For $k = 2$, we note that the $p$-completion of the $p$-units in
  $K_n$ is topologically generated by $\zeta_n$ and those $p$-units fixed by
  complex conjugation.  But this means only even eigenspaces
  and the eigenspace of $\mc{P}_{L/K}^{(2)}$ generated by
  $$
    (\lambda_{2r-1},\zeta)_{p^{\infty},K_{\infty}}^{(2)}
    \otimes \zeta^{\otimes -2}
  $$
  can be nontrivial.  Galois acts on the latter element by the
  $(3-4r)$th power of the cyclotomic character.  On the other hand,
  this element generates a subgroup of $\epsilon_{1-r}A_K \cong \zp(1-r)$
  by Lemma \ref{consequence}.  Since $3-4r = 1-r$ if and only if
  $3r = 2 \in T$, we are done.
\end{proof}

\begin{remark}
  The condition of $\lambda$-invariant $1$ implicitly
  assumed in Proposition \ref{larger} is also known to hold for all
  nontrivial eigenspaces of $\mc{A}_K$ for $p <$ 12,000,000
  \cite{bcems}.  It is not expected to hold in general.
\end{remark}

Venjakob \cite[Section 8]{venjakob} and Coates \cite[Section
3]{coates} asked if $\mc{Z}_L$ is always pseudo-null as a
$\zp[[\mc{G}]]$-module (see \cite[Definition 3.1]{venjakob-str}
for the definition).  By \cite[Proposition 6.5]{venjakob}, the
$\zp[[\mc{G}]]$-module $\mc{Z}_L$ is pseudo-null if and only if it
(or, equivalently, $\mc{A}_{L}$) is torsion as a $\zp[[G]]$-module
(since $\mc{Z}_L$ is finitely generated as a $\zp[[G]]$-module).
This would be a natural extension of a conjecture of Greenberg's
\cite[Conjecture 3.5]{greenberg} to this noncommutative case. As
is demonstrated in a forthcoming paper with Hachimori \cite{hs},
this is not so if
we take $L_{\infty}$ to be a CM-field and the
$(-1)$-eigenspace $\mc{A}_K^-$ of $\mc{A}_K$ under complex
conjugation has $\Z_p$-rank at least $2$.  In fact, a result of
considerably greater generality is proven.  We give another proof
of the above-mentioned result (slightly weakened for simplicity)
here.

\begin{proposition} \label{pseudo}
    Let $L_{\infty}$ be a CM $\zp$-extension of $K_{\infty}$ which is
    Galois over $K$, unramified outside $p$, and totally ramified at $p$.
    Then $\mc{A}_L^-$ is not pseudo-null if the
    $\zp$-rank of $\mc{A}_K^-$ is at least $2$.
\end{proposition}

\begin{proof}
    Note that $G$ is fixed under complex conjugation, since $L$ is
    CM.  It is easily checked that
    \begin{equation*}
        \rank_{\zp[[G]]} \mc{A}_L^- = \lim_{k \to \infty} \rank_{\zp} 
        \mr{gr}_G^k \mc{A}_L^-.
    \end{equation*}
    By Theorem \ref{mainthm}, this yields
    \begin{equation} \label{limitrank}
        \rank_{\zp[[G]]} \mc{A}_L^- = \rank_{\zp} \mc{A}_K^- -
        \lim_{k \to \infty} \rank_{\zp} (\mc{P}_{L/K}^{(k)})^-,
    \end{equation}
    since $L/K$ is totally ramified at $p$.
    Now $\mc{U}_K^-$ is simply the Tate module of $p$-power roots
    of unity, and the Kummer generator $a$
    of $L_{\infty}$ lies in
    the $(-1)$-eigenspace of the pro-$p$ completion of
    $K_{\infty}^{\times}$.  Hence $(\mc{P}_{L/K}^{(k)})^-$ is
    generated over $\zp$ by $(a,\zeta)^{(k)}_{p^{\infty},K_{\infty},S} \otimes_{\zp}
    \zeta^{\otimes -k}$.  By \eqref{limitrank}, we are done.
\end{proof}

Note that $\mc{A}_K^-$ has $\zp$-rank (at least) 2 for $p = 157,
353, 379, \ldots$. The CM condition implies that $L_{\infty}$ is
{\em not} generated by a sequence of $p$-power roots of $p$-units
over $K_{\infty}$. In the situation that $L_{\infty}/K_{\infty}$
is generated by cyclotomic $p$-units, it is still an open question
whether $\mc{Z}_L$ is always pseudo-null or not. We phrase it as a
conjecture (see \cite{hs} for a related question).

\begin{conjecture} \label{unitconj}
    Let $L_{\infty}$ be a $\zp$-extension of $K_{\infty} = \Q(\mu_{p^{\infty}})$
    unramified outside $p$ and generated over $K_{\infty}$ by $p$-power roots of
    $p$-units in $K_{\infty}$.  Then $\mc{Z}_L$ is pseudo-null.
\end{conjecture}

Proposition \ref{p37} provides examples in which $\mc{Z}_{L}$ is
indeed pseudo-null.  In \cite{me-eis}, we establish a sufficient
and computable condition for the nontriviality of the cup product
pairing which holds at least for $p < 1000$, in particular
furnishing a much larger class of examples of this form.  Beyond
pseudo-nullity, it would be interesting to know if it ever occurs
that $I_G^k\mc{Z}_L$ is nontrivial for all $k$ when $L_{\infty}$
arises from roots of $p$-units.

We end by remarking on a possible counterexample to Conjecture
\ref{unitconj}.

\begin{remark}
In \cite[Proposition 3.1]{wingberg} and \cite[Exemple 3.2]{nqd},
it was claimed that the Galois group of the maximal pro-$p$
abelian unramified extension of $K_{\infty}$ (in which all primes
above $p$ split completely) is free pro-$p$ if Vandiver's
conjecture holds at $p$. Under Vandiver's conjecture, any
unramified $\zp$-extension $L_{\infty}/K_{\infty}$ is generated by
$p$-power roots of cyclotomic $p$-units. Assuming the freeness
result, it is easily seen that $\mc{Z}_L$ is not pseudo-null
whenever the $\zp$-rank of $\mc{A}_K$ is at least $2$ and
Vandiver's conjecture holds at $p$. (We thank Manabu Ozaki for
pointing out this example.) Unfortunately, however, the proofs of
the above-mentioned result both have subtle but vital
mistakes.  For example, the Galois group of
the maximal unramified pro-$157$ extension of $K_{\infty}$ for
$p=157$ is actually isomorphic to $\Z_{157} \oplus \Z_{157}$ (note that $157$
divides $B_{62}$ and $B_{110}$). This follows from Proposition
\ref{smallclassgp} and the nontriviality of the cup product
$\langle \lambda_{47,1},\lambda_{15,1} \rangle_{157,62}$
(for this nontriviality, see \cite{mcs,me-eis}).
\end{remark}

\renewcommand{\baselinestretch}{1}

\footnotesize \noindent
Department of Mathematics and Statistics\\
McMaster University\\
1280 Main Street West\\
Hamilton, Ontario L8S 4K1\\ 
Canada\\
email address: {\tt sharifi@math.mcmaster.ca}
\end{document}